\documentclass[10pt]{amsart}

\usepackage{euscript,graphicx,epstopdf,amscd,amsgen,amsfonts,amssymb,latexsym,
amsmath,amsthm,graphicx,mathrsfs,times,color,overpic}

\newtheorem{theorem}{Theorem}[section]
\newtheorem*{theorem*}{Theorem}
\newtheorem{lemma}[theorem]{Lemma}
\newtheorem{scholium}[theorem]{Scholium}
\newtheorem{corollary}[theorem]{Corollary}
\newtheorem{claim}[theorem]{Claim}
\newtheorem{proposition}[theorem]{Proposition}

\numberwithin{equation}{section}

\def\cV{\mathcal V}
\def\scA{\mathscr A}

\def\cU{\mathcal U}
\def\var{{\rm var}}

\def\s{{\rm s}}
\def\ss{{\rm ss}}
\def\uu{{\rm uu}}

\def\ss{{\rm ss}}
\def\u{{\rm u}}

\newtheorem{maintheorem}{Theorem}

\theoremstyle{definition}
\newtheorem{definition}[theorem]{Definition}

\newtheorem*{example*}{Example}
\newtheorem{remark}[theorem]{Remark}
\newtheorem*{remark*}{Remark}

\newcommand{\eqdef}{\stackrel{\scriptscriptstyle\rm def}{=}}

\def\s{{\rm s}}
\def\u{{\rm u}}
\def\uu{{\rm uu}}

\DeclareMathOperator{\interior}{int}

\DeclareMathOperator{\card}{card}
\DeclareMathOperator{\proj}{proj}

\def\bA{\mathbb{A}}

\def\bN{\mathbb{N}}

\def\bZ{\mathbb{Z}}
\def\bR{\mathbb{R}}

\def\cC{\mathscr{C}}

\def\cF{\mathcal{F}}

\def\cO{\EuScript{O}}

\def\cP{\EuScript{P}}

\def\cU{\EuScript{U}}
\def\cV{\EuScript{V}}

\def\cM{\EuScript{M}}

\DeclareMathSymbol{\varnothing}{\mathord}{AMSb}{"3F}
\author[K. Gelfert]{Katrin Gelfert}\address{Institute of Mathematics, Federal University of Rio de Janeiro, Av. Athos da Silveira Ramos 149, Cidade Universit\'aria - Ilha do Fund\~ao, Rio de Janeiro 21945-909, RJ, Brazil}\email{gelfert@im.ufrj.br}
\author[F. Riquelme]{Felipe Riquelme}\address{IMA, Pontificia Universidad Cat\'olica de Valpara\'iso, Blanco Viel 596, Valpara\'iso, Chile.}\email{felipe.riquelme@pucv.cl}
\thanks{F.R. was supported by FONDECYT Iniciaci\'on N$^o$11190461 and CAPES (Brazil). K.G. was supported [in part] by CAPES -- Finance Code 001 and CNPq (Brazil).}
\begin{document}

\title[Exceptional sets]{Exceptional sets for geodesic flows\\ of noncompact manifolds}

\begin{abstract} For a geodesic flow on a negatively curved Riemannian manifold $M$ and some subset $A\subset T^1M$, we study the limit $A$-exceptional set, that is the set of points whose $\omega$-limit do not intersect $A$. We show that if the topological $\ast$-entropy of $A$ is smaller than the topological entropy of the geodesic flow, then the limit $A$-exceptional set has full topological entropy. Some consequences are stated for limit exceptional sets of invariant compact subsets and proper submanifolds.
\end{abstract}

\begin{thanks}{}\end{thanks}
\maketitle

\section{Introduction}


Exceptional sets have been largely studied over the last years, mostly trying to identify for which dynamical systems small sets have large exceptional sets. Naturally, there are plenty of ways to quantify the ``size'' of a set, but these studies have been focused mainly in two notions: Hausdorff dimension, which involves only the topology of a set, and topological entropy, which mostly considers the dynamic contribution of a set on the system. 

In the discrete-time case Dolgopyat \cite{Dol:97} proved that sets with small entropy in subshifts of finite type have exceptional sets with full entropy. Analogous statements hold for piecewise expanding maps of the interval and for Anosov diffeomorphisms of the two-dimensional torus, but in terms of Hausdorff dimension (see \cite{Dol:97}). Campos and Gelfert \cite{CaGe:16,CaGe:19} concluded the same type of results in terms of both, topological entropy and Hausdorff dimension, for nonuniformly expanding/hyperbolic maps. To the best of our knowledge the first result involving exceptional sets for continuous-time dynamical systems comes from \cite{Kle:98} where the author proves that the exceptional set of a compact (proper) submanifold on a negatively curved Riemannian manifold having constant negative curvature has full Hausdorff dimension relative to the geodesic flow. A similar result can be found in \cite{KleWei:13} for partially hyperbolic flows on homogeneous spaces. Most of these results were motivated by the work of Jarnik-Besicovitch \cite{Ja:29} on the estimation of the size of badly approximable real numbers and its equivalence with the study of the set of bounded geodesic orbits on the modular surface. Indeed, the set of bounded orbits can be informally thought as the exceptional set of $\infty$. Dani \cite{Dan:86} proved that the set of bounded orbits by the action of a one-parameter subgroup on $G/\Gamma$ has full Hausdorff dimension, where $G$ is a connected semisimple Lie group of real rank 1 and $\Gamma$ is a lattice (see also \cite{KleMar:96} for analogous results on homogeneous dynamics). 

The aim of this work is to study the topological entropy of exceptional sets for geodesic flows on arbitrary noncompact negatively curved Riemannian manifolds, this includes those with variable sectional curvatures and dimension $\ge 2$.

To state our first main result, let us introduce some notation. Consider a semi-group $F=(f^t)_{t\in \bA}$, $\bA=\bR_{\ge0}$ or $\bA=\bN\cup\{0\}$, acting on some complete metric space $X$. Denote by $\cO^+_{F}(x)=\{f^t(x)\colon t\in\bA\}$ the (\emph{forward}) \emph{semi-orbit} of $x\in X$ by $F$. Denote by $\omega_F(x)$ the (\emph{forward}) \emph{$\omega$-limit set} of a point $x\in X$, that is, the set of limit points of $\cO^+_{F}(x)$. Given a set $A\subset X$, denote by
\[
	I_F(A)
	\eqdef \{x\colon \omega_{F}(x)\cap A=\emptyset\}
\]
the (\emph{forward}) \emph{limit $A$-exceptional set}  (with respect to $F$). Notice that $I_F(A)$ is $F$-invariant, that is $f^t(I_F(A))=I_F(A)$ for every $t\in\bA$. In the case when $F|_Y\colon Y\to Y$ is the restricted flow on some $F$-invariant set $Y$, then we use the notation $I_{F|Y}(A)$. 

Given a set $A\subset X$, we denote by $h^\ast_{\rm top}(F,A)$ the \emph{topological $\ast$-entropy} of $F$ on $A$,
\[
	h_{\rm top}^\ast(F,A)
	\eqdef \sup h_\mu(F),
\] 
where the supremum is taken over all Borel probability measure which are weak$\ast$ limits of empirical measures distributed along the forward orbit of a point in $A$, and we let $h_{\rm top}^\ast(F,A)\eqdef0$ in case no such measure exists. By $h_{\rm top}(F,A)$ we denote the \emph{topological entropy} of $F$ on $A$.  We  recall the precise definition of entropy in Section \ref{sec:entropy}.

 The following is our first main result. 

\begin{maintheorem}\label{the:1}
	Consider the geodesic flow $G=(g^t)_{t\in\bR}$ on the unit tangent bundle $T^1M$ of a $n$-dimensional complete Riemannian manifold $M$ of negatively pinched sectional curvatures $-b^2\le \kappa\leq -1$ for some $b\geq 1$. For every Borel set $A\subset T^1M$ satisfying $h_{\rm top}^\ast(G,A)<h_{\rm top}(G,T^1M)$ it holds
\[
	h_{\rm top}^\ast(G,I_{G}(A))
	= h_{\rm top}(G,T^1M).
\]	
\end{maintheorem}

Let us provide some (nontrivial) example of Borel sets satisfying the hypothesis of Theorem \ref{the:1}. Let $D$ be the set of \emph{divergent in average geodesic orbits}, that is, orbits which spend a Birkhoff average time 0 on every compact set. Then $h_{\rm top}^\ast(G,D) =0$ since there are no limit measures for empirical measures distributed along divergent in average orbits. In particular, $h_{\rm top}^\ast(G,I_{G}(D))
	= h_{\rm top}(G,T^1M)$. Note that the set $D$ is not so small since by \cite{RiVe:22} it has positive Hausdorff dimension provided $M$ is not convex-cocompact. 
	
	Let us state now some consequences of Theorem \ref{the:1}. Recall that the \emph{non-wandering set} $\Omega$ of the geodesic flow $G$ is the set of unit vectors $v\in T^1M$ such that for every neighborhood $U$ of $v$ there exist $s,t>0$ such that $g^s(v),g^{-t}(v)\in U$. 

\begin{corollary}\label{cor:1}
	Assume the hypotheses of Theorem \ref{the:1} and additionally that the derivatives of the sectional curvatures are uniformly bounded. 
	If the non-wandering set of the geodesic flow is non-compact, then for every compact and $G$-invariant set $K\subset T^1M$,
\[
	h^\ast_{\rm top}(G,I_{G}(K)) 
	= h_{\rm top}(G).
\]
\end{corollary}

For the following consequence, consider a proper Riemannian submanifold $N\subset M$. Then $N$ is isometric to the quotient $N=\widetilde{N}/\Gamma'$, where $\widetilde{N}$ is a Riemannian submanifold of the universal covering of $M$, $\widetilde{M}$, and $\Gamma'$ is a subgroup of $\Gamma$ acting on $\widetilde{N}$. A discrete group $H\subset {\rm Isom}^+(\widetilde{M})$ is \emph{divergent} if its Poincar\'e series  
\[
	\mathcal{P}_H(s)
	=\sum_{h\in H} e^{-sd(o,h o)}
\]
diverges at $s=\delta_H$, where $\delta_H$ is its \emph{critical exponent}, 
\[
	\delta_H
	=\limsup_{R\to\infty}\frac{1}{R}\log\card\{h\in H \colon d(o,h o)\leq R\}.
\]
Denote by $\Omega_N$ the non-wandering set associated to the geodesic flow over $T^1N$.

\begin{corollary} \label{cor:2}
		Assume the hypotheses of Theorem \ref{the:1} and additionally that the derivatives of the sectional curvatures are uniformly bounded. 
	Let $N\subset M$ be a proper Riemannian submanifold as above. Assume that $\Gamma'$ is divergent and $\Omega_N\neq \Omega$. Then
\[
h^\ast_{\rm top}(G,I_{G}(T^1N)) = h_{\rm top}(G).
\]
\end{corollary} 

In particular, Corollary \ref{cor:2} applies if $M$ is compact with constant negative curvature $-1$ and $N$ is a proper submanifold of $M$. Indeed, in this case $\Gamma'$ is divergent and $\Omega_N=T^1N\neq  T^1M=\Omega$. We provide the proofs of the above results in Section \ref{sec:proofs}.

The proof of Theorem \ref{the:1} will rely on the consideration of appropriate suspension flows, namely, we construct basic sets (that is, sets which are compact, $G$-invariant, and locally maximal such that $G|_B$ is topologically transitive) having almost full entropy, where the flow is conjugated to a suspension flow over a subshift of finite type (see Proposition \ref{pro:entropybasicset}). In particular, the following is an immediate consequence.

\begin{scholium}\label{sch:1}
	Assume the hypotheses of Theorem \ref{the:1}. Denote by $B$ the set of vectors $v\in T^1M$ whose orbit $\cO_G(v)\eqdef\{g^t(v)\colon t\in\bR\}$ is bounded. Then
\[
	\sup_{K\subset T^1M\text{ basic}}h_{\rm top}(G,K)
	= h_{\rm top}(G,B)
	= h^\ast_{\rm top}(G,B)
	= h_{\rm top}(G,T^1M).
\]	 
\end{scholium}

We finally state the corresponding auxiliary key result, which is of independent interest. 

\begin{maintheorem}\label{thepro:susflo}
	Consider some continuous function $\tau\colon\Sigma\eqdef\{1,\ldots,N\}^\bZ\to (0,\alpha]$ and let $F\colon\Sigma(\sigma,\tau)\times\bR\to\Sigma(\sigma,\tau)$ be the suspension flow of the shift map $\sigma\colon\Sigma\to\Sigma$ under $\tau$. Then for every Borel set $A\subset\Sigma(\sigma,\tau)$ satisfying $h^\ast_{\rm top}(F,A)<h_{\rm top}(F)$ it holds 
\[
        h^\ast_{\rm top}(F,I_{F}(A))
	= h_{\rm top}(F).
\]	
\end{maintheorem}

This paper is organized as follows. In Section \ref{sec:entropy} we recall different concepts of entropy and some fundamental properties. Section \ref{sec:geodynneg} introduces the necessary elements from geometry. In Section \ref{sec:Katokshift} we prove an approximation result by ergodic measures on the symbolic space.
In Section \ref{sec:Katokflow}, we study techniques of ``ergodic approximation by basic sets'' of an ergodic probability measure for the geodesic flow of a complete manifold (see Proposition \ref{pro:entropybasicset}). These approximation techniques obtaining flow-invariant basic sets will be implemented in Section \ref{sec:susflo} in which we study exceptional sets for suspension flows. Theorem \ref{thepro:susflo} is proven in Section \ref{sec:susflo}. Theorem \ref{the:1}, Corollaries \ref{cor:1} and \ref{cor:2}, and Scholium \ref{sch:1} are proven in Section \ref{sec:proofs}.

\section{Entropy}\label{sec:entropy}

We study topological entropy on arbitrary (not necessarily invariant or compact) subsets where the ambient space is not necessarily compact. Hence, in this section we recall some essential definitions and results. 

Given a continuous flow $F=(f^t)_{t\in\bR}$ on a metric space $(X,d)$, a Borel probability measure $\mu$ is \emph{$F$-invariant} if $(f^t)_\ast\mu=\mu$ for every $t\in\bR$. Denote by $\cM(F)$ the set of all $F$-invariant Borel probability measures. A measure $\mu\in\cM(F)$ is \emph{ergodic} if for every $F$-invariant set $B\subset X$  either $\mu(B)=0$ or $\mu(B)=1$. Denote by $\cM_{\rm erg}(F)$ the subset of ergodic measures. 
Given $t\in\bR$, we denote by $\cM(f^t)$ the set of $f^t$-invariant Borel probability measures and by $\cM_{\rm erg}(f^t)$ the subset of $f^t$-ergodic ones.

\begin{remark}\label{rem:ergodic}
	Note that for every $\mu\in\cM_{\rm erg}(F)$ for every time $t>0$, except for a countable set of $t$-values values,  it holds that $\mu\in\cM_{\rm erg}(f^t)$ \cite{PugShu:71}.
\end{remark}

\subsection{Metric entropy}
Given a measure $\mu\in\cM(F)$, we denote by $h_\mu(F)$ its \emph{metric entropy} (with respect to the flow $F$) and by $h_\mu(f^t)$ its entropy (with respect to the map $f^t$). By Abramov's formula \cite{Abr:59}, for every $t\in\bR$ it holds 
\begin{equation}\label{eq:Abramov}
	\lvert t\rvert \, h_\mu(f^1)=h_\mu(f^t)
\end{equation}
and we let $h_\mu(F)\eqdef h_\mu(f^1)$.

\subsection{Topological entropy}
Let $(X,d)$ be a metric space and $f\colon X\to X$ a continuous map. For $n\in\bN$ and $\varepsilon>0$ define 
\begin{equation}\label{eq:Bowbal}
	B_n^d(x,\varepsilon)
	\eqdef \{y\in X\colon d(f^k(x),f^k(y))<\varepsilon\text{ for all }k\in\{0,\ldots,n-1\}\}
\end{equation}
and call it \emph{$(n,\varepsilon)$-dynamical ball}.
The above is simply the ball of radius $\varepsilon$ centered at $x$ with respect to the metric $d_n$ defined by
\[
	d_n(x,y)
	\eqdef \max_{k\in\{0,\ldots,n-1\}}d(f^k(x),f^k(y)).
\]
A set $S\subset X$ is \emph{$(n,r)$-separated} if $d_n(x,y)\ge r$ for every pair of points $x,y\in S$, $x\ne y$. As $f$ is continuous, $B_n^d(x,\varepsilon)$ is open. Given a compact set $K\subset X$, denote by $N^d(n,\varepsilon,K)$  the minimal cardinality of a cover of $K$ by $(n,\varepsilon)$-dynamical balls. Then
\[
	h^d_{\rm top}(f)
	\eqdef \sup_K\lim_{\varepsilon\to0}\limsup_{n\to\infty}\frac1n\log N^d(n,\varepsilon,K),
\]
where the supremum is taken over all compact subsets $K\subset M$, is the \emph{topological entropy} of $f$ (with respect to $d$). 
If the metric $d$ is clear from the context, we will drop the superscript ${}^d$.

Metric entropy and topological entropy are linked by the variational principle (see \cite{Din:70}  for $X$ compact and \cite{HanKit:95} for the general case). 
First note that 
\[
	\sup_{\mu\in\cM(f)} h_\mu(f)
	\le h^d_{\rm top}(f).
\]
A change of metric may result in a change of entropy. Assuming that $X$ is a locally compact metrizable topological space, the following \emph{variational principle} holds
\begin{equation}\label{eq:vp}
	\sup_{\mu\in\cM(f)} h_\mu(f)
	= h_{\rm top}(f)
	\eqdef \inf_d h^d_{\rm top}(f),
\end{equation}
where the infimum is taken over all metrics $d$ on $X$ that generate the same topology. One calls $h_{\rm top}(f)$ the \emph{topological entropy} of $f$. Note that, by ergodic decomposition, in the supremum in \eqref{eq:vp} it is enough to consider the ergodic measures only. In the case when $(X,d)$ is a compact metric space then $h_{\rm top}(f)=h^d_{\rm top}(f)$.

By \eqref{eq:Abramov} together with the variational principle \eqref{eq:vp},  for a continuous flow $F=(f^t)_{t\in\bR}$ on a locally compact metrizable topological space, for every $t\in\bR$ it holds
\[
	h_{\rm top}(f^t)
	=\lvert t\rvert\, h_{\rm top}(f^1),
\]
and we define the \emph{topological entropy} of $F$ by
\begin{equation}\label{eq:timetmap}
	h_{\rm top}(F)
	\eqdef h_{\rm top}(f^1).
\end{equation}
The following variational principle is a consequence of the above facts
\begin{equation}\label{eq:vpflow}
	h_{\rm top}(F)
	= \sup_{\mu\in\cM(f^1)}h_\mu(F)
	=  \sup_{\mu\in\cM(F)}h_\mu(F).	
\end{equation}

\subsection{Topological entropy on Borel subsets}\label{sec:entmeasub}

Let $f\colon X\to X$ be a continuous map on a (not necessarily compact) metric space $(X,d)$. Given $x\in X$ and $n\in\bN$, consider the probability measure
\begin{equation}\label{eq:deltadef}
	\delta_{x,n}
	\eqdef \frac1n\sum_{k=0}^{n-1}\delta_{f^k(x)}.
\end{equation}
Denote by 
\[
	\cV(f,x)
	\eqdef \{\mu\in\cM(f)\colon \delta_{x,n_k}\to\mu\text{ for some }n_k\to\infty\}
\]
the set of limit points relative to the weak$\ast$ topology. Recall that a sequence of measures $(\mu_n)$ converge weak$\ast$ to $\mu$ if and only if for every bounded continuous function $\varphi:X\to \bR$, we have $\lim \int \varphi d\mu_n=\int \varphi d\mu$. Note that, since $X$ is not necessarily compact, the set $\cV(f,x)$ could be empty. If the map $f$ is clear from the context then we drop it in the notation. Following \cite{Tho:11}, given a nonempty Borel set $Z\subset X$, denote 
\begin{equation}\label{notationMZ}
	\cM^Z(f)
	\eqdef \{\mu\colon \mu\in\cV(f,x)\text{ for some }x\in Z\}
\end{equation}
and define by
\begin{equation}\label{def:entstar}
	h^\ast_{\rm top}(f,Z)
	\eqdef 
	\begin{cases}
		\sup\big\{h_\mu(f)\colon \mu\in\cM^Z(f)\}&\text{ if }\cM^Z(f)\ne\emptyset,\\
		0&\text{ otherwise}.
	\end{cases}
\end{equation}
the \emph{topological entropy} of $f$ on $Z$.

We will implement the above definition in the noncompact case only to state Theorem \ref{the:1}. For the remaining arguments, we will consider appropriate compact and invariant subset $Y\subset X$ and study the entropy of $f|_Y$ on $Z\cap Y$. Note that together with \cite[Theorem 6.1]{Tho:11}, in this case it holds
\[
	h^\ast_{\rm top}(f|_Y,Z\cap Y)
	= h^\ast_{\rm top}(f,Z\cap Y)
	\le h^\ast_{\rm top}(f,Z).
\]

\subsection{Topological entropy on subsets}\label{sec:entcomsub}

We also recall the definition of entropy according to Bowen \cite{Bow:73} for a continuous map $f\colon X\to X$ on a general topological space. Given an open cover $\scA$ of $X$ and a subset $Z\subset X$, denote by $n_\scA(Z)$ the smallest nonnegative integer $n$ such that $f^n(Z)$ is not contained in an element of $\scA$; if $f^n(Z)$ is contained in an element of $\scA$ for all integers $n\ge0$ then let $n_\scA(Z)=\infty$. Define
\[
	m_\scA(Z,s)
	\eqdef \lim_{r\to0}\inf_\cU\sum_{U\in\cU}e^{-sn_\scA(U)},
\]
where the infimum is taken over all countable open covers $\cU$ of $Z$ such that $n_\scA(U)>1/r$ for all $U\in\cU$. The \emph{topological entropy} of $f$ on $Z$  is defined by
\[
	h_{\rm top}(f,Z)
	\eqdef \sup_\scA h_\scA(f,Z),
	\quad\text{ where }\quad
	h_\scA(f,Z)
	\eqdef \inf\{s\colon m_\scA(Z,s)=0\}.
\]
If $Z=X$ is compact then the latter equals the topological entropy $h_{\rm top}(f)$ defined in \eqref{eq:vp}.

We collect below some elementary properties relating all previous notions of entropy.

\begin{lemma}\label{lem:propTho}
Assume that $(X,d)$ is a metric space and $f\colon X\to X$ a continuous map. 
\begin{enumerate}
\item[(i)] 
Given a nonempty Borel set $Z\subset X$, for every $n\in\bN$ it holds
\[
	h^\ast_{\rm top}(f^n,Z)
	= n\, h^\ast_{\rm top}(f,Z).
\]
\item[(ii)] 
Given nonempty Borel sets $Y,Z\subset X$ satisfying $Y\subset Z$, it holds
\[
	h^\ast_{\rm top}(f,Y)
	\le h^\ast_{\rm top}(f,Z).
\]
\end{enumerate}
\end{lemma}

\begin{proof}
	Recall  $h_\mu(f^n)=n h_\mu(f)$ and note that $\cV(f^n,x)\subset\cV(f,x)$. This together implies $h^\ast_{\rm top}(f^n,Z) \le n h^\ast_{\rm top}(f,Z)$. On the other hand,  if $\mu\in\cV(f,x)$ then for every $n\in\bN$ it holds $\nu\eqdef\frac1n(\mu+f_\ast\mu+\ldots+f^{n-1}_\ast\mu)\in\cV(f^n,x)$ and $h_\nu(f^n)=h_\mu(f^n)=n\,h_\mu(f)$. This proves item (i). Item (ii) is an immediate consequence of the definition.
\end{proof}

\begin{lemma}\label{lem:propTho2}
Assume that $(X,d)$ is a compact metric space, $f\colon X\to X$ is a continuous map, and $Z\subset X$ is a Borel set.  
\begin{enumerate}
\item[(i)] If there are a compact metric space $(Y,\rho)$, a continuous map $g\colon Y\to Y$, and a continuous surjective map $p\colon X\to Y$ satisfying $p\circ f=g\circ p$, then 
\[
	h_{\rm top}(g,p(Z))
	\le h_{\rm top}(f,Z).
\]	
\item[(ii)] If the map $p$ in item (i) is a homeomorphism, then 
\[
	h^\ast_{\rm top}(f,Z)= h^\ast_{\rm top}(g,p(Z)) 
	\quad\text{ and }\quad
	h_{\rm top}(f,Z)= h_{\rm top}(g,p(Z)) .
\]	
\item[(iii)] It holds
\[
	h_{\rm top}(f,Z)
	\le h_{\rm top}^\ast(f,Z).
\]
\item[(iv)] It holds
\[
	h_{\rm top}(f,X)= h_{\rm top}^\ast(f,X).
\]
\end{enumerate}
\end{lemma}

\begin{proof}
	The proof of item (i) is straightforward. Item (ii)  is \cite[Theorem 3.2]{Tho:11}, item (iii) is \cite[Theorem 4.3]{Tho:11}, and item (iv) is \cite[Theorem 3.3]{Tho:11}.
\end{proof}

\begin{remark} We stress the fact that inequality (iii) in Lemma \ref{lem:propTho2} is strict in general. Take for instance an ergodic measure $\mu$ having positive entropy and $x\in X$ a $\mu$-generic point. Then for $Z=\{x\}$ we have $h_{\rm top}(f,Z)=0$ whereas $h_{\rm top}^\ast(f,Z)=h_\mu(f)>0$.
\end{remark}

In analogy to the above, for a continuous flow $F=(f^t)_{t\in\bR}$ on a metric space $X$, we define the \emph{topological entropy} of $F$ on a nonempty Borel set $Z\subset X$ by
\[
	h_{\rm top}^\ast(F,Z)
	\eqdef h_{\rm top}^\ast(f^1,Z).
\]

\section{Geometry and dynamics in negative curvature}\label{sec:geodynneg}

Let $\widetilde{M}$ be a complete simply connected Riemannian manifold with sectional curvatures bounded from above by $-1$. The \emph{boundary at infinity} $\partial_\infty\widetilde{M}$ of $\widetilde{M}$ is the set of asymptotic geodesic rays on $\widetilde{M}$. In particular, the set $\widetilde{M}\cup\partial_\infty\widetilde{M}$ is compact endowed with the cone topology and homeomorphic to the closed unit ball. We stress the fact that every isometry of $\widetilde{M}$ extends to a homeomorphism of $\widetilde{M}\cup\partial_\infty\widetilde{M}$.

\subsection{Hopf parametrization} Fix once for all a point $o\in\widetilde{M}$. For any boundary point $\xi\in\partial_\infty\widetilde{M}$, let $r_\xi\colon [0,+\infty)\to\widetilde{M}$ be the arc-parametrization of the geodesic ray with origin $o$ and extremity at infinity $\xi$. The \emph{Busemann cocycle} of $\widetilde{M}$ is the map $\beta\colon\widetilde{M}\times\widetilde{M}\times\partial_\infty\widetilde{M}\to\bR$ defined by
\[
(x,y,\xi)\mapsto\beta_\xi(x,y)\eqdef\lim_{t\to+\infty}d(x,r_\xi(t))-d(y,r_\xi(y)).
\]
This limit exists by the bounds on the sectional curvatures. It is independent of $o$. 


Let $T^1\widetilde{M}$ be the unit tangent bundle of $\widetilde{M}$ and $\pi\colon T^1\widetilde{M}\to\widetilde{M}$ be the natural projection. For every $v\in T^1\widetilde{M}$, let $v^-$ and $v^+$ be the two extremities at infinity of the geodesic line defined by $v$. Let $\partial^2_\infty\widetilde{M}$ be the subset of $\partial_\infty\widetilde{M}\times\partial_\infty\widetilde{M}$ consisting of distinct points at infinity. The \emph{Hopf parametrization} of $T^1\widetilde{M}$ is the identification of $v\in T^1M$ with the triplet $(v^-,v^+,s)\in\partial^2_\infty\widetilde{M}\times\bR$, where $s=\beta_{v^+}(o,\pi(v))$. This map is a homeomorphism. We  use the notation $v=(v^-,v^+,s)$ whenever we mention a vector $v\in T^1\widetilde{M}$ in these coordinates. 

\subsection{Distances on (sub)manifolds} 

Let $v\colon(-\infty,+\infty)\to\widetilde M$ be the parametrization of the oriented geodesic ray defined by $v\in T^1\widetilde{M}$ such that $\pi(v)=v(0)$ and $\pi v\colon t\mapsto \pi(v(t))$ is arc-parametrized. We endow $T^1\widetilde{M}$ with the distance $d$ defined for all $v,v'\in T^1\widetilde{M}$ as follows
\[
	d(v,v')
	\eqdef\frac{1}{\sqrt{\pi}}\int d(\pi(v(t)),\pi(v'(t)))e^{-t^2}dt.
\]

\begin{remark} 
The distance $d$ is H\"older-equivalent to the distance induced by the Sasaki metric on $T^1\widetilde{M}$ \cite[Lemma 2.3]{PauPolSch:15}. 
\end{remark}

We denote by $G\eqdef (g^t)_{t\in\bR}$ the geodesic flow on $T^1\widetilde{M}$ and $\iota\colon T^1\widetilde{M}\to T^1\widetilde{M}$ the flip map $v\mapsto -v$. 
Note that the geodesic flow acts by translation in the third coordinate, namely $g^t(v^-,v^+,s)=(v^-,v^+,s+t)$, and the flip map sends $(v^-,v^+,s)$ into $(v^+,v^-,-s)$. Both, the geodesic flow and the flip map commute with every isometry of $\widetilde{M}$. Moreover $\iota\circ g^t = g^{-t}\circ\iota$ for every $t\in\bR$.

\begin{lemma}\label{lem:closegeodesics} 
Let $\varepsilon\in(0,1]$ and $s\geq 2$. There exists a constant $C=C(\varepsilon)>0$ such that, for any two geodesic lines $(\xi,\eta)$ and $(\xi',\eta')$ on $\widetilde{M}$ staying $\varepsilon$-close for time of length at least $2s$, their Hausdorff distance is less than $Ce^{-s}$. 
\end{lemma}

\begin{proof}
Let $p\in (\xi,\eta)$ and $p'\in (\xi',\eta')$ be such that $d(p,p')$ realizes the Hausdorff distance. Let $r,r'\colon\bR\to \widetilde{M}$ be the arc length parametrizations of $(\xi,\eta)$ and $(\xi',\eta')$, respectively, such that $r(0)=p$ and $r'(0)=p'$. By hypothesis, $d(r(t),r'(t))\leq\varepsilon$ for all $-s\leq t\leq s$. Set $x=r(-s)$, $y=r(s)$, $x'=r'(-s)$, and $y'=r'(s)$. Let $[z,w]$ be the shortest arc between $[x,x']$ and $[y,y']$ with $z\in [x,x']$. Let $q\in\widetilde{M}$ be the midpoint of $[z,w]$ and $\ell=\frac{1}{2}d(z,w)$. Compare also Figure \ref{figure1}. Note that $q$ is also the midpoint of the segment defined by $p$ and $p'$. 

\begin{figure}
\begin{overpic}[scale=.4
  ]{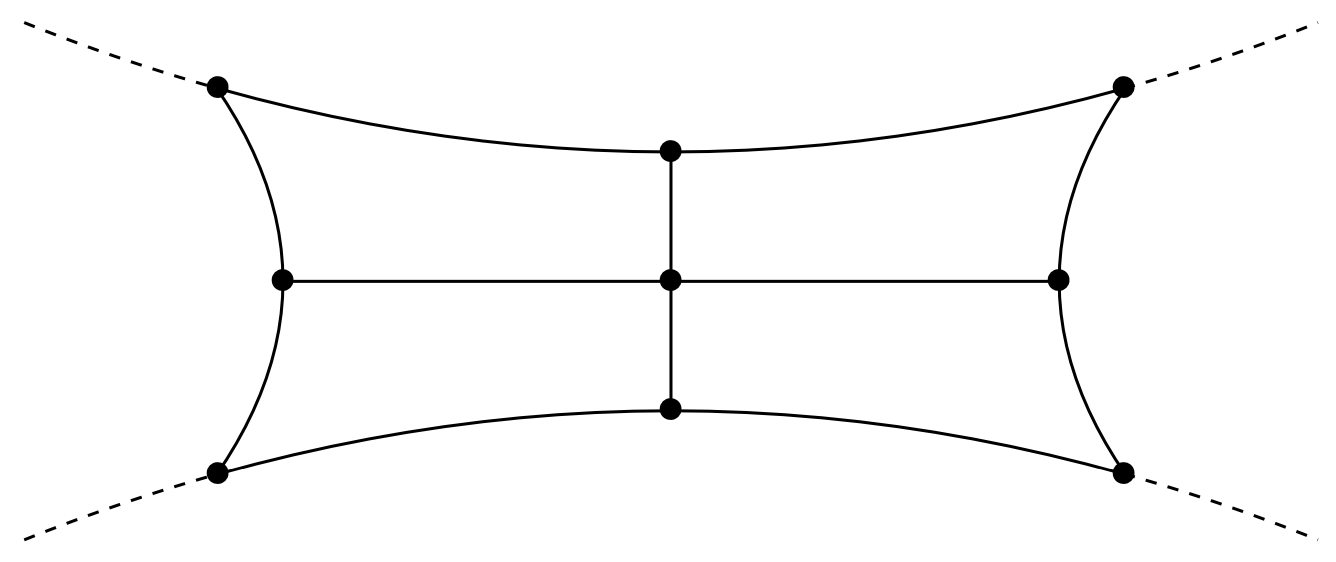}
      	\put(50,34){$p$}	
      	\put(50,7){$p'$}
      	\put(65,33){$s$}	
      	\put(35,33){$s$}	
      	\put(65,8){$s$}	
      	\put(35,8){$s$}	
      	\put(65,18){$\ell$}	
      	\put(35,18){$\ell$}	
      	\put(18,18){$z$}	
      	\put(80,18){$w$}	
      	\put(16,3){$x'$}	
      	\put(85,38){$y$}	
      	\put(17,37){$x$}	
      	\put(84,2){$y'$}	
      	\put(-5,20){\small $d(x,x')\le\varepsilon$}	
      	\put(90,20){\small $d(y,y')\le\varepsilon$}	
      	\put(-2,40){$\xi$}	
      	\put(-2,0){$\xi'$}	
      	\put(102,0){$\eta'$}	
      	\put(102,40){$\eta$}	
      	\put(51,18){$q$}	
\end{overpic}
\caption{Proof of Lemma \ref{lem:closegeodesics}}
\label{figure1}
\end{figure}

By comparison, the distance $d(p,q)$ is less than the distance $\tau$ in the hyperbolic plane $\mathbb{H}$ between the midpoint $\bar{q}$ of a segment $[\bar{z},\bar{w}]$ of length $2l$ to a geodesic segment $[\bar{x},\bar{y}]\subset\mathbb{H}$ where the angle at $\bar{z}$ between $[\bar{x},\bar{z}]$ and $[\bar{z},\bar{w}]$ and the angle at $\bar{w}$ between $[\bar{z},\bar{w}]$ and $[\bar{w},\bar{y}]$ are exactly $\pi/2$. By \cite[Theorem 7.17.1 (i)]{Bea:83}, it holds $\sinh(\tau)\sinh(\ell)\leq 1$.
 Hence
\[
	d(p,p')
	=2d(p,q)
	\leq 2\tau
	\leq2\sinh(\tau)
	= 2\sinh(\ell)^{-1}.
\]
By the triangle inequality, we also have $2\ell=d(z,w)\geq d(x,y)-2\varepsilon=2(s-\varepsilon)$, hence $\ell\geq s-\varepsilon$. By the inequalities above, we therefore have
\[
d(p,p')\leq 2\sinh(s-\varepsilon)^{-1}.
\]
On the other hand, the hyperbolic sine function verifies $\sinh(\alpha)\geq e^{\alpha}/4$ for every $\alpha\geq 1$. Since $s-\varepsilon\geq 1$, we get
\[
d(p,p')\leq 8e^{-(s-\varepsilon)}=8e^{\varepsilon}e^{-s},
\]
which ends the proof of the lemma.
\end{proof}

The \emph{strong stable} and \emph{strong unstable manifolds} of $v$ are defined by
\[\begin{split}
	\widetilde{W}^{\ss}(v)
	&\eqdef \{w\in T^1\widetilde{M}\colon d(v(t),w(t))\to 0 \quad \text{as} \quad t\to+\infty\},\\
	\widetilde{W}^{\uu}(v)
	&\eqdef \{w\in T^1\widetilde{M}\colon d(v(t),w(t))\to 0 \quad \text{as} \quad t\to-\infty\},
\end{split}\]
respectively.
These submanifolds are smooth leaves of continuous foliations, namely the \emph{strong stable} and \emph{strong unstable foliations} $\widetilde{W}^{\ss}$ and $\widetilde{W}^{\uu}$, which are invariant under the isometry group of $\widetilde{M}$ and the geodesic flow. Note that, in Hopf coordinates, for every $v=(v^-,v^+,s)\in T^1\widetilde{M}$, 
\[\begin{split}
	\widetilde{W}^{\ss}(v)
	&=\{(w^-,w^+,r)\in T^1\widetilde{M} \colon w^+=v^+, \ r= s\},\\
	\widetilde{W}^{\uu}(v)
	&=\{(w^-,w^+,r)\in T^1\widetilde{M} \colon w^-=v^-, \ r= s\}.
\end{split}\]

\begin{definition} 
The \emph{Hamenst\"adt distances} are defined, for every $u\in T^1\widetilde{M}$, as 
\[\begin{split}
	d^\u_u(v,v')
	&\eqdef\lim_{t\to+\infty} e^{\frac{1}{2}d(\pi(v(t)),\pi(v'(t)))-t}
	\quad	\text{ for all $v,v'\in \widetilde{W}^{\uu}(u)$}.\\
	d^\s_u(w,w')
	&\eqdef\lim_{t\to+\infty}e^{\frac{1}{2}d(\pi(w(-t)),\pi(w'(-t)))-t}
	\quad	\text{ for all $w,w'\in \widetilde{W}^\ss(u)$}.
\end{split}\] 
Note that both limits exists and defines a distance inducing the original topology on $\widetilde{W}^\uu(u)$ and $\widetilde{W}^\ss(u)$ respectively (see \cite{Ham:89} for further details). 
\end{definition}

Note that 
\[
	\iota(\widetilde{W}^{\uu}(u))=\widetilde{W}^{\ss}(\iota(u))
\]
It is straightforward to check that $d^\u_u(v,v')=d^\s_{\iota u}(\iota v,\iota v')$ for all $v,v'\in \widetilde{W}^{\uu}(u)$. 

 The following lemma appears in \cite[Lemma 2.4]{PauPolSch:15}, following \cite{ParkPau:14}.

\begin{lemma}\label{lem:hamen:uu} 
There is $c>0$ so that for every $u\in T^1\widetilde{M}$, $v,v'\in \widetilde{W}^{\uu}(u)$, and $w,w'\in\widetilde W^\ss(u)$,
\[\begin{split}
	\max\left\{\frac{1}{c}d(v,v'),d(\pi(v),\pi(v'))\right\}
	&\leq d^\u_u(v,v')\leq e^{\frac{1}{2}d(\pi(v),\pi(v'))},\\
	\max\left\{\frac{1}{c}d(w,w'),d(\pi(w),\pi(w'))\right\}
	&\leq d^\s_u(w,w')\leq e^{\frac{1}{2}d(\pi(w),\pi(w'))}.
\end{split}\]
\end{lemma}

%

A remarkable property of the Hamenst\"adt distances is the uniform expansion/contraction under the action of the geodesic flow. 

\begin{lemma}\label{lem:unifexp}
For every $u\in T^1\widetilde{M}$,  $s\in\bR$, $v,v'\in \widetilde{W}^{\uu}(u)$, and $w,w'\in \widetilde{W}^{\ss}(u)$, 
\[\begin{split}
	d^{\u}_{g^s(u)}(g^s(v),g^s(v'))
	&=e^s d^{\u}_u(v,v')\\
	d^{\s}_{g^s(u)}(g^s(w),g^s(w'))
	&=e^{-s} d^{\s}_u(w,w').
\end{split}\]
\end{lemma}

\subsection{Quotients} Let $\Gamma$ be a discrete, non-elementary subgroup of isometries of $\widetilde{M}$ and $M\eqdef\widetilde{M}/\Gamma$ the quotient space. Even if $M$ is not a manifold (since $\Gamma$ is not necessarily torsion-free), we denote by $T^1M$ the quotient $T^1\widetilde{M}/\Gamma$. Since the geodesic flow and the flip map commute with every isometry of $\widetilde{M}$, both descend respectively to maps on $T^1M$ that we still denote as $g^t\colon T^1M\to T^1M$ and $\iota\colon T^1M\to T^1M$. For $v\in T^1M$ we define the \emph{strong stable} and \emph{strong unstable manifolds} $W^{\ss}(v)$ and $W^{uu}(v)$ as
\[\begin{split}
	W^{\ss}(v)
	\eqdef\{w\in T^1M\colon d(g^t(v),g^t(w))\to 0 \quad \text{as} \quad t\to+\infty\},\\
	W^{\uu}(v)
	\eqdef\{w\in T^1M\colon d(g^t(v),g^t(w))\to 0 \quad \text{as} \quad t\to-\infty\},
\end{split}\]
respectively.
Let $p^1_\Gamma\colon T^1\widetilde{M}\to T^1M$ and $p_\Gamma\colon \widetilde{M}\to M$ be the corresponding quotient maps. Then for any $\tilde{v}\in T^1\widetilde{M}$ and $v=p^1_\Gamma(\tilde{v})$, it holds $p^1_\Gamma(\widetilde{W}^{\ss}(\tilde{v}))\subset W^{\ss}(v)$ and $p^1_\Gamma(\widetilde{W}^{\uu}(\tilde{v}))\subset W^{\uu}(v)$. For small $\varepsilon>0$, we also define the \emph{local strong stable manifold} $W^{\ss}_\varepsilon(v)$ (resp., \emph{local strong unstable manifold} $W^{\uu}_\varepsilon(v)$) at $v\in T^1M$ as the connected component of $W^{\ss}(v)\cap B(v,\varepsilon)$ (resp. $W^{\uu}(v)$) containing $v$. We can similarly define local strong stable/unstable manifolds on $T^1\widetilde{M}$, and denote them by $\widetilde{W}^{\ss}_\varepsilon(\tilde{v})$ and $\widetilde{W}^{\uu}_\varepsilon(\tilde{v})$. Observe that for $\varepsilon$ small enough (depending on the injectivity radius at $\pi(v)$), it holds
\[
p^1_\Gamma(\widetilde{W}^{\ss}_\varepsilon(\tilde{v}))=W^{\ss}_\varepsilon(v) \quad \text{and} \quad p^1_\Gamma(\widetilde{W}^{\uu}_\varepsilon(\tilde{v}))=W^{\uu}_\varepsilon(v),
\]
so strong stable/unstable manifolds can be locally studied using coordinates in $T^1\widetilde{M}$. We also remark that Hamenst\"adt distances are inherited on
 local strong stable/unstable manifolds on $T^1M$.

\subsection{Local product structure} 

The geodesic flow on negatively curved manifolds verifies the following local product structure.

\begin{figure}
\begin{overpic}[scale=.4
  ]{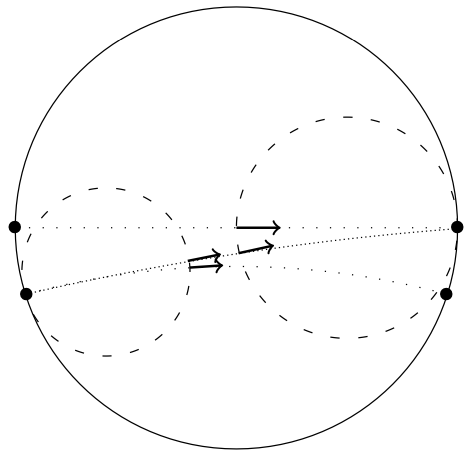}
      	\put(44,80){$\widetilde M$}	
      	\put(89,18){$\partial_\infty\tilde M$}	
      	\put(100,48){\small$\widetilde u^+$}	
      	\put(98,32){\small$\widetilde v^+$}	
      	\put(-3,48){\small$\widetilde u^-$}	
      	\put(-2,32){\small$\widetilde v^-$}	
      	\put(44,34){\small$\widetilde v$}	
      	\put(55,50){\small$\widetilde u$}	
      	\put(58,40){\small$\widetilde w$}	
\end{overpic}
\caption{Local product structure}
\label{figure}
\end{figure}
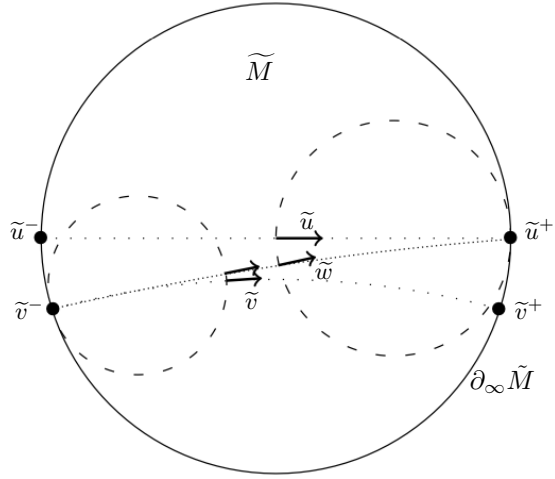

\begin{proposition}[Local product structure]\label{prop:lps}  Every $u\in T^1M$ admits a neighbourhood $V$ which satisfies the following. For every $\varepsilon>0$ there exists $\beta>0$ such that for all $v\in V$ with $d(u,v)<\beta$, there exists a unique vector $w\in T^1M$ and a real number $|t|<\varepsilon$, so that
\[
w\eqdef{[u,v]}\in W^{\ss}_\varepsilon(u)\cap W^{\uu}_\varepsilon(g^t v)
\]
\end{proposition}

\begin{proof}
Let $u\in T^1M$, $\tilde{u}\in T^1\widetilde{M}$ any lift, and let
\[
	V_{\tilde{u}}
	\eqdef \{\tilde{v}\in T^1\widetilde{M}\colon \tilde{v}^-\neq \tilde{u}^+\}.
\]
The neighborhood $V_{\tilde{u}}$ of $\tilde{u}$ is open and dense in $T^1\widetilde{M}$. Moreover, any $\tilde{v}\in V_{\tilde{u}}$ has the following property: there exists $t\in\bR$, $\tilde{v}_{\ss}\in \widetilde{W}^{\ss}(\tilde{u})$ and $\tilde{v}_{\uu}\in \widetilde{W}^{\uu}(\tilde{v})$ such that $g^t(\tilde{v}_{\uu})=\tilde{v}_\ss$, that is $\widetilde{W}^{\ss}(\tilde{u})\cap \widetilde{W}^{\uu}(g^t\tilde{v})$ is non-empty with only one element $\tilde{w}\eqdef\tilde{v}_\ss$. Since strong stable/unstable leaves are continuous, as well as the geodesic flow, for any $\varepsilon>0$ there exists $\beta>0$ such that $d(\tilde{u},\tilde{v})<\beta$ implies $\tilde{v}\in V_{\tilde{u}}$, $d(\tilde{u},\tilde{w})<\varepsilon$, $d(g^t(\tilde{v}),\tilde{w})<\varepsilon$ and $|t|<\varepsilon$. Finally, the desired property is verified in $T^1M$ for $w=p^1_\Gamma(\tilde{w})$ whenever $\varepsilon$ is chosen less than the injectivity radius at $u$.  Compare also Figure \ref{figure}. 
\end{proof}



\section{Approximating ergodic measures by subshifts of finite type}\label{sec:Katokshift}

In this section, we follow techniques introduced in \cite[Supplement S.5]{KatHas:95} to construct some compact invariant subset which ergodically (weak$\ast$ and in entropy) approximates a given ergodic measure. We invoke this approach in the case of a subshift of finite type. It will be useful when proving Theorem \ref{thepro:susflo} in Section \ref{sec:susflo}.

Given $N\in\bN$, consider the shift space $\Sigma\eqdef\{1,\ldots,N\}^\bZ$ equipped with the metric $d_\Sigma(\underline i,\underline j)\eqdef 2^{-\inf\{\lvert k\rvert\colon i_k\ne j_k\}}$ and the left shift $\sigma\colon\Sigma\to\Sigma$. A \emph{subshift of finite type} (\emph{SFT}) is a subset  $\Xi\subset\Sigma$ for which there is a matrix $A=(a_{k\ell})_{k,\ell=1}^N$, $a_{k\ell}\in\{0,1\}$,  so that $\Xi=\{\underline i\in\Sigma\colon a_{i_ni_{n+1}}=1\text{ for all }n\in\bZ\}$. Note that any SFT is compact and $\sigma$-invariant.  For every $\underline i=(i_0i_1\ldots)\in\Sigma$ and $n\in\bN$, denote the \emph{$n$th level cylinder} containing $\underline i$ by
\[
	[i_0\ldots i_{n-1}]
	\eqdef \{\underline j\in\Sigma\colon j_k=i_k\text{ for all }k=0,\ldots,n-1\}.
\]

	Let $\pi^+\colon\Sigma\to\Sigma^+$ be the natural projection defined by $\pi^+(\underline i)\eqdef\underline i^+\eqdef(i_0i_1\ldots)$. Consider also the space $\Sigma^+\eqdef\{\underline i^+\colon \underline i\in\Sigma\}$ of one-sided sequences together with the left shift $\sigma^+\colon\Sigma^+\to\Sigma^+$. A SFT for $\sigma^+$ is analogously defined. Analogously, consider the $n$th cylinder of $\underline i^+$, which we denote by $[i_0\ldots i_{n-1}]^+$. Let $\nu^+\eqdef\pi^+_\ast\nu$ and note that $\nu^+$ is $\sigma^+$-ergodic. Moreover,  $h_{\nu^+}(\sigma^+)=h_{\nu}(\sigma)$.

\begin{proposition}\label{pro:SFT}
	For every continuous function $\tau\colon\Sigma\to (0,\alpha]$, $\nu\in\cM_{\rm erg}(\sigma)$, and $\varepsilon>0$ there exists a subshift of finite type $\Xi\subset\Sigma$ such that  
\[
	\lvert h_{\rm top}(\sigma,\Xi)-h_\nu(\sigma)\rvert
	<\varepsilon
	\quad\text{ and }\quad
	\Big\lvert \int\tau\,d\lambda-\int\tau\,d\nu\Big\rvert
	<\varepsilon
	\quad\text{for every $\lambda\in\cM(\sigma|_\Xi)$}.
\]
\end{proposition}

Before proving the above result, for $n\in\bN$ and $\psi\colon\Sigma\to\bR$ let
\[
	\var_n\psi
	\eqdef \sup\{\lvert \psi(\underline j)-\psi(\underline i)\rvert\colon \underline j,\underline i\in\Sigma, j_k=i_k\text{ for all }\lvert k\rvert\le n\}
\]
and define the set of functions
\[
	\cF_\Sigma
	\eqdef \{\psi\colon\Sigma\to\bR\colon
	 \text{there exist } b>0,\alpha\in(0,1)\text{ so that }
	\var_n\psi\le b\alpha^n\text{ for all }n\ge0\}.
\] 

\begin{proof}[Proof of Proposition \ref{pro:SFT}] 
	Let $\varepsilon>0$. 	
Let $\psi\in\cF_\Sigma$ such that $\sup\lvert\psi-\tau\rvert<\varepsilon/3$. Indeed, for example any function which is piecewise constant on sufficiently high-level cylinders has this property.

	By \cite[1.6 Lemma]{Bow:08} there is $\phi\in\cF_\Sigma$ which is cohomologous to $\psi$ and satisfies $\phi(\underline k)=\psi(\underline i)$ for every $\underline k\in\Sigma$ such that $k_n=i_n$ for all $n\ge0$. Hence the function $\psi^+\colon\Sigma^+\to\bR$, $\psi^+(\underline i^+)\eqdef \phi(\underline i)$ for any $\underline i\in[i_0i_1\ldots]$ is well-defined and continuous. As $\phi\in\cF_\Sigma$, there are $b>0$ and $\alpha\in(0,1)$ such that for all $n\ge0$ it holds $\var_n\phi\le b\alpha^n$. Hence, from the definition of $\psi^+$, we get $\var_n\psi^+\le b\alpha^n$.
		
	By the Brin-Katok theorem \cite{BriKat:83}, for $\nu^+$-almost every $\underline i^+\in\Sigma^+$ it holds
\[
	\lim_{n\to\infty}-\frac1n\log\nu^+([i_0\ldots i_{n-1}]^+)
	= h_{\nu^+}(\sigma^+)
	= h_{\nu}(\sigma)
	\eqdef h.
\]	
	By the Birkhoff ergodic theorem, for $\nu^+$-almost every $\underline i^+\in\Sigma^+$ it holds
\[
	\lim_{n\to\infty}\frac1nS_n\psi^+(\underline i^+)
	= \int\psi^+\,d\nu^+,
	\quad\text{ where }
	S_n\psi^+\eqdef \psi^++\psi^+\circ \sigma^++\ldots+\psi^+\circ(\sigma^+)^{n-1}.
\]

	Fix $\kappa\in(0,1)$.
By the Egorov theorem, there is a set $\Lambda\subset\Sigma^+$ satisfying $\nu^+(\Lambda)>1-\kappa$ and $N\in\bN$ such that for every $n\ge N$ and $\underline i^+\in\Lambda$ it holds
\begin{equation}\label{eq:cylin}
	e^{-n(h-\varepsilon)}
	\le \nu^+([i_0\ldots i_{n-1}]^+)
	\le e^{-n(h+\varepsilon)}
\end{equation}
and 
\[
	\left\lvert \int\psi^+\,d\nu^+
		-\frac1nS_n\psi^+(\underline i^+)\right\rvert
	\le \frac\varepsilon2.
\]
Assume that $N$ was chosen large enough such that $b\alpha^N\le\varepsilon/2$. Hence for every $\underline i^+\in\Lambda$ and $n\ge N$,
\[
	\max_{\underline j^+\in[i_0\ldots i_{n-1}i_n\ldots i_{2n-1}]^+}
		\left\lvert S_n\psi^+(\underline j^+)
			-S_n\psi^+(\underline i^+)\right\rvert
	\le nb\alpha^n
	\le n\frac\varepsilon2	.
\]
Fix now some $n\ge N$.
Choose any sequence $\underline i^{+,1}\in \Lambda$. Let $A_0\eqdef \Lambda$. Inductively, for $k\in\bN$ choose $\underline i^{+,k}\in A_{k-1}$ and let $A_k\eqdef A_{k-1}\setminus[i^k_0\ldots i^k_{n-1}]^+$. By \eqref{eq:cylin} there is $M\in\bN$,
\[
	e^{n(h-\varepsilon)}
	\le M
	\le e^{n(h+\varepsilon)},
\] 
such that for every $k>M$ the set $A_k$ is empty.  By definition, the cylinders $[i^{k}_0\ldots i^{k}_{n-1}]^+$, $k=1,\ldots,M$, are pairwise disjoint and 
\[
 	(\sigma^+)^n\big([i^{k}_0\ldots i^{k}_{n-1}]^+\big)
	=\Sigma^+.
\]

The (one-sided) infinite concatenation of any combination of finite sequences from the family $\{(i^{k}_0\ldots i^{k}_{n-1})\}_{k=1}^M$ defines a (one-sided) SFT $\Xi^+\subset\Sigma^+$. By the above
\[
	h-\varepsilon
	\le \frac1n \log \card M
	= h_{\rm top}(\sigma^+,\Xi^+) 	
	\le h+\varepsilon .
\]

For any ergodic measure $\lambda^+\in\cM(\sigma^+|_{\Xi^+})$, for $\lambda^+$-almost every $\underline i^+$ it holds
\[
	\int\psi^+\,d\lambda^+
	= \lim_{k\to\infty}\frac{1}{kn}S_{kn}\psi^+(\underline i^+).
\] 
Hence 
\[
	\Big\lvert \int\psi^+\,d\lambda^+-\int\psi^+\,d\nu^+\Big\rvert
	<\varepsilon.
\]
Let $\sigma|_\Xi\colon\Xi\to\Xi$ be the natural extension of the SFT $\Xi^+$ and note that it is SFT (with respect to $\sigma$) and satisfies $h_{\rm top}(\sigma,\Xi)=h_{\rm top}(\sigma^+,\Xi^+)$. The above implies that for every 
ergodic measure $\lambda\in\cM(\sigma|_{\Xi})$, 
\[\begin{split}
	\varepsilon
	&>\Big\lvert \int\psi\,d\lambda-\int\psi\,d\nu\Big\rvert\\
	&\ge \Big\lvert \int\tau\,d\lambda-\int\tau\,d\nu\Big\rvert
	- \Big\lvert \int\tau\,d\nu-\int\psi\,d\nu\Big\rvert
	-\Big\lvert \int\psi\,d\lambda-\int\tau\,d\lambda\Big\rvert\\
	&\ge \Big\lvert \int\tau\,d\lambda-\int\tau\,d\nu\Big\rvert
		- 2\frac\varepsilon3.
\end{split}\]
This proves the proposition. 
\end{proof}

\section{Approximating basic sets in a geodesic flow}\label{sec:Katokflow}

In this section we continue the principle idea in Section \ref{sec:Katokshift} to ``ergodically approximate'' an ergodic measure. We deal now with the geodesic flow of a complete Riemannian manifold of negative curvature and obtain an approximation in terms of entropy by basic sets.
For that first recall that a compact $G$-invariant set $B\subset T^1M$ is \emph{locally maximal} if there exists a neighborhood $U$ of $B$ such that
\[
	B
	=\bigcap_{t\in\bR}g^t(\overline U).
\]
Recall that $G|_B$ is \emph{topologically transitive} if there is $x\in B$ such that $\overline{\cO^+_G(x)}=B$. We say that a compact $G$-invariant set is \emph{basic} if it is hyperbolic, locally maximal, and $G|_B$ is topologically transitive. 
A set $B\subset T^1M$ has \emph{a local product structure} if  for every $u,v\in B$ satisfying $d(u,v)<\beta$ for $\beta>0$ sufficiently small the point $[u,v]$ is again in $B$ (where $[\cdot,\cdot]$ was defined in Proposition \ref{prop:lps}).
Recall that a compact $G$-invariant hyperbolic set $B$ is locally maximal if and only if it has a local product structure (see, for example, \cite[Theorem 6.2.7]{FisHas:19}).

Note that the construction of compact invariant sets with certain ergodic properties in this context can also be found in \cite[Section 6.2]{PauPolSch:15} while proving a variational principle for the topological pressure, following ideas from \cite[Section 4]{OP:04}. 

The following is our main result of this section.

\begin{proposition}\label{pro:entropybasicset}
	Consider the geodesic flow $G=(g^t)_{t\in\bR}$ on the unit tangent bundle $T^1M$ of a $n$-dimensional complete Riemannian manifold $M$ of negatively pinched sectional curvatures $-b^2\le \kappa\leq -1$ for some $b\geq 1$. 
	Then for every measure $\mu\in\cM_{\rm erg}(G)$ with positive entropy and $\varepsilon>0$ there exists a basic set $B\subset T^1M$ of topological dimension $1$ satisfying
\[
	h_{\rm top}(G,B)
	\ge h_\mu(G)-\varepsilon.
\]	 
\end{proposition}

We postpone the proof of this proposition to Section \ref{ssec:proofprop}. We first prepare the necessary tools. In Section \ref{ssec:prorec} we first consider sections which induce proper rectangles, obtained from the local product structure. In Section \ref{ssec:prorecK} we cover a given compact set by proper rectangles and study a return map to local cross sections.

\subsection{Sections and proper rectangles}\label{ssec:prorec}

By the local product structure stated in Proposition \ref{prop:lps}, for every $\varepsilon=\varepsilon_{\rm lps}>0$ sufficiently small, there is $\beta_{\rm lps}>0$ such that for every $u,v\in T^1M$ satisfying $d(u,v)\le\beta_{\rm lps}$  the point
\begin{equation}\label{eq:locpro}
	w
	\eqdef [u,v]
	= W^\ss_\varepsilon(u)\cap W^\uu_\varepsilon(g^t(v))
	\quad\text{ for some }\lvert t\rvert\le\varepsilon.
\end{equation}
is well defined (this intersection contains just one point). 

Let us now transfer this local product structure to local cross sections of the flow.
Given an interval $I\subset\bR$ and $Z\subset T^1M$, denote $g^I(Z)\eqdef \bigcup_{t\in I}g^t(Z)$.
A \emph{section of size $\varepsilon$} is a closed set $D\subset T^1M$ such that  $(x,t)\mapsto g^t(x)$ is a homeomorphism between $D\times[-\varepsilon,\varepsilon]$ and its image $g^{[-\varepsilon,\varepsilon]}(D)$. In the following, we will assume that $D$ is a sufficiently small closed codimension-one smooth disk. Every section $D$ has associated a well-defined projection map $\proj_D\colon g^{[-\varepsilon,\varepsilon]}(D)\to D$ given by $\proj_D(g^t(u))=u$. Note that the domain of $\proj_D$ contains a nonempty open subset of $T^1M$.

If $D$ is a section, then for a closed set $R\subset D$ satisfying $d(R,\partial D)>0$ and having diameter sufficiently small relative to $d(R,\partial D)$, we define 
\[
	[\cdot,\cdot]_D\colon R\times R\to D,\quad 
	[u,v]_D\eqdef\proj_D([u,v]).
\]
We say that $R$ is a \emph{rectangle} in $D$ if $[R,R]_D\subset R$; in this case let $[\cdot,\cdot]_R\eqdef[\cdot,\cdot]_D|_{R\times R}$. A rectangle $R$ in $D$ is \emph{proper} if $R=\overline{\interior R}$ in the internal topology of the disk $D$.

For $D$ a section of size $\varepsilon$ and $R$ a proper rectangle in $D$, consider the flow-box 
\[
	\widehat R
	\eqdef g^{[-\varepsilon/2,\varepsilon/2]}(R).
\]
This way, for every $u\in \widehat R$ there is $\tau=\tau(R,u)\in[-\varepsilon/2,\varepsilon/2]$ so that $g^{\tau}(u)\in R\subset D$. Given $w\in R$, let
\[\begin{split}
	W^\s(w,R)
	&\eqdef \{[w,v]_R\colon v\in R\}
	= R\cap\proj_D\big(g^{[-\varepsilon,\varepsilon]}(D)
		\cap W^\ss_\varepsilon(w)\big),\\
	W^\u(w,R)
	&\eqdef \{[v,w]_R\colon v\in R\}
	= R\cap\proj_D\big(g^{[-\varepsilon,\varepsilon]}(D)
		\cap W^\uu_\varepsilon(w)\big).
\end{split}\]
This way, $R$ is a direct product of the sets $W^\s(w,R)$ and $W^\u(w,R)$ in terms of the homeo\-morphism $(u,v)\mapsto[u,v]_R$. 

We state the following lemma without proof (see, for example, the arguments for constructing geometric rectangles in \cite[Section 4.1]{ConLafTho:20}).

\begin{lemma}
	Given any $v\in T^1M$, there exist a proper rectangle $R$ of arbitrarily small diameter such that $v\in \interior R$.	
\end{lemma}

\subsection{Sections and proper rectangles to cover a compact set}\label{ssec:prorecK}

For the statement of the following lemma, recall the constant $c>0$ from Lemma \ref{lem:hamen:uu}. 

\begin{lemma}\label{lem:rec1}
	For every compact set $K\subset T^1M$ and $\alpha>0$, for every $\varepsilon_{\rm sec}\in(0,\alpha/(3c))$ sufficiently small, there are a finite collection of sections $\{D_i\}$ of size $\varepsilon_i\in(0,\varepsilon_{\rm sec})$ and a number $\beta_{\rm sec}>0$ such that for every $\beta\in(0,\beta_{\rm sec})$, there are points $v_i\in K$ such that
\[
	K
	\subset \bigcup_i B(v_i,\frac\beta2)
\]	
and that for every $i$ there is a proper rectangle $R_i$ in $D_i$ such that $B(v_i,\beta)\subset g^{(-\varepsilon_i,\varepsilon_i)}(R_i)$. 

Moreover, there is $M\in\bN$ such that for every index $i$, for every $m\ge M$ the following holds. The map
\begin{equation}\label{def:Tim}\begin{split}
	&T_{i,m}\colon\widehat D_i^m\to \widehat D_i^m,
	\quad\text{ where }\quad
	\widehat D_i^m\eqdef \proj_{D_i}\big(\widehat R_i \cap g^m(\widehat R_i)\big)\\
	&T_{i,m}(u')
	\eqdef g^{m+\tau(R_i,g^m(u))}(u')
	\in R_i,
	\quad\text{ where }\quad
	u'=g^{\tau(R_i,u)}(u)
\end{split}\end{equation}
is continuous and injective and
\begin{equation}\label{eq:estTim}
	T_{i,m}(u')
	= g^t(u')
	\quad\text{ for some }\quad
	t\in(m-\varepsilon_i,m+\varepsilon_i).
\end{equation}
Moreover, denoting by $\cC_{i,m}(u)$ the connected component of $\widehat D_i^m$ which contains $u\in D_i$, every such component is a proper rectangle in $D_i$.  
Furthermore, if $v,w\in B(v_i,\beta/2)\cap g^{-m}(B(v_i,\beta/2))$ are $(m,\alpha)$-separated, then 
\begin{equation}\label{eq:disjoint}
	\cC_{i,m}(\proj_{D_i}(v))\cap\cC_{i,m}(\proj_{D_i}(w))
	=\emptyset.
\end{equation}	
\end{lemma}

\begin{proof}
Let $\varepsilon_{\rm lps}\in(0,\alpha/(3c))$ be sufficiently small and choose $\beta_{\rm lps}>0$ accordingly to the local product structure. Fix
\[
	0<\varepsilon_{\rm sec}
	< \frac{1}{3}\min\{c,\frac1c\}\cdot\min\{\varepsilon_{\rm lps},\alpha\}.
\]		

	By compactness of $K$ and continuity of the local product structure, for every $u\in K$ there exists a smooth section $D(u)$ of size $\varepsilon_u\in(0,\varepsilon_{\rm sec})$ containing $u$ and a proper rectangle $R(u)$ in $D(u)$. Note that $U(u)\eqdef g^{(-\varepsilon_u/2,\varepsilon_u/2)}(\interior R(u))$ defines an open neighborhood of $u$ in $T^1M$. By compactness, there exists a finite cover $\{U(u_i)\colon u_i\in K\}$ of $K$. Let $L_0>0$ be a Lebesgue number for this open cover. 
	
Let 
\[
	\beta\in(0,\beta_{\rm sec}),
	\quad\text{ where }\quad
	\beta_{\rm sec}
	\eqdef\min\{c,\frac1c\}\cdot\min\{\frac12L_0,\beta_{\rm lps}\}.
\]	
Fix some finite open cover of $K$, 
\[
	K
	\subset \bigcup_iB(v_i,\frac\beta2).
\]	
By the above, for every index $i$ there is some $u_i$ so that 
\[
	B(v_i,\beta)
	\subset g^{(-\varepsilon_i/2,\varepsilon_i/2)}(\interior R_i)
	\subset\widehat R_i,
\]	 
where $\varepsilon_i\eqdef\varepsilon(u_i)$ and $R_i\eqdef R(u_i)$ is a proper rectangle in the section $D_i\eqdef D(u_i)$ of size $\varepsilon_i$. This proves the first assertion of the lemma.

Choose now $M\in\bN$ sufficiently large such that 
\begin{equation}\label{hyp:betaM}
	e^{-M}e^{\varepsilon_{\rm lps}}\varepsilon_{\rm lps}<\frac1c\frac\beta2.
\end{equation}	 

\begin{claim}\label{cla:prorec}
	For every index $i$, for every $m\ge M$ and $u\in B(v_i,\beta/2)\cap g^{-m}(B(v_i,\beta/2))$,
\[\begin{split}
	W^\s\big(\proj_{D_i}(u),R_i\big)
	&\subset \proj_{D_i}\big(\widehat R_i\cap g^{-m}(\widehat R_i)\big),\\
	W^\u\big(\proj_{D_i}(g^m(u)),R_i\big)
	&\subset \proj_{D_i}\big(g^m(\widehat R_i)\cap \widehat R_i\big).
\end{split}\]
\end{claim}

\begin{proof}
As $u\in B(v_i,\beta/2)\subset\widehat R_i$,  $u'=g^s(u) \eqdef\proj_{D_i}(u)$ for $s\eqdef \tau(R_i,u)$ satisfying $\lvert s\rvert\le\varepsilon_i/2$. Given $w'\in W^\s(u',R_i)$. Hence, by definition, $w'\in R_i\subset\proj_{D_i}(\widehat R_i)$. 

On the other hand, also by definition, $w'=g^t(w)=\proj_{D_i}(w)$ for some $w\in W^\ss_{\varepsilon_i}(u')$ and $t\eqdef\tau(R_i,w)$ satisfying $\lvert t\rvert\le\varepsilon_i$. Therefore, $w''\eqdef g^{-s}(w)\in g^{-s}(W^\ss(u'))= W^\ss(g^{-s}(u'))=W^\ss(u)$. By Lemma \ref{lem:unifexp},
\[\begin{split}
	d^\s_{g^m(u)}(g^m(w''),g^m(u))
	&= e^{-m}d^\s_{u}(w'',u)
	= e^{-m}d^\s_u(g^{-s}(w),g^{-s}(u'))\\
	&= e^{-m}e^sd^\s_{u'}(w,u')
	\le e^{-m}e^{\varepsilon_i}\varepsilon_i\\
	\text{\small by \eqref{hyp:betaM}}\quad
	&< \frac1c\frac\beta2.
\end{split}\]
By Lemma \ref{lem:hamen:uu},
\[
	d(g^m(w''),g^m(u))
	\le cd^\s_{g^m(u)}(g^m(w''),g^m(u))
	\le c\frac1c\frac\beta2
	= \frac\beta2.
\]
Hence, our hypothesis $g^m(u)\in B(v_i,\beta/2)$ implies that $g^m(w'')\in B(v_i,\beta)\subset\widehat R_i$. In particular, $w'=\proj_{D_i}(w'')\in g^{-m}(\widehat R_i)$. This proves the first inclusion.

 The other inclusion is analogous.
\end{proof}

Claim \ref{cla:prorec} hence implies that every set $\cC_{i,m}(\proj_{D_i}(u))$ is a proper rectangle in $D_i$.

For every index $i$, for every $m\ge M$ and $u\in \widehat R_i\cap g^{-m}(\widehat R_i)$ the image $g^m(u)$ is in the flow-box $\widehat R_i$, and hence the orbit of $u$ passes through the rectangle $R_i$. Hence, its associated return  to $R_i$, \eqref{def:Tim}, is well defined. Note that it is not necessarily the first return. 

\begin{claim}\label{cl:Timcontinuous}
	The map $T_{i,m}$ defined in \eqref{def:Tim} is continuous and injective.
\end{claim}	
	
\begin{proof}
	To show continuity, note that $u\mapsto \tau(R_i,u)$ is continuous on $\widehat R_i$. Hence, continuity of $T_{i,m}$ follows from the continuity of the flow $g$.

	By contradiction, assume that $g^{m+\tau(R_i,u)}(u)=g^{m+\tau(R_i,v)}(v)\eqdef w\in R_i$ for $u,v\in R_i$, $u\ne v$. Hence, letting $u'=g^m(u)$ and $v'=g^m(v)$, it holds $v'=g^\delta(u')$ for some $\delta\in[-\varepsilon_i,\varepsilon_i]$, $\delta\ne0$. Without loss of generality, we can assume that $\delta>0$. Then, because $u',v'$ are both contained in the flow-box $\widehat R_i$, it follows that the segment $g^{[0,\delta]}(u')$ is completely contained in this flow-box. Since on one hand $u,v\in \widehat R_i$ and on the other hand $u,v\in g^{-m}(g^{[0,\delta]}(u')) =g^{[-m,m+\delta]}(u')$, it follows that $u,v$ are on the same orbit. But this contradicts the fact that $u,v\in D_i$ and the fact that $D_i$ is a section.
\end{proof}

%
	For every $u'\in\widehat D_i^m$, it holds
\[
	T_{i,m}(u')
	= T_{i,m}(g^{\tau(R_i,u)}(u))
	= g^{m+\tau(R_i,g^m(u))-\tau(i,u)}(u'),
\]	
where $\tau(R_i,u),\tau(R_i,g^m(u))\in[-\varepsilon_i/2,\varepsilon_i/2]$, proving
\eqref{eq:estTim}. 

Let us finally show \eqref{eq:disjoint}.

\begin{claim}\label{cla:difcomp}
	For every index $i$, for every $m\ge M$ and pair of points $v,w\in B(v_i,\beta/2)\cap g^{-m}(B(v_i,\beta/2))$ which are $(m,\alpha)$-separated, it holds
\[
	\cC_{i,m}(\proj_{D_i}(v))\cap\cC_{i,m}(\proj_{D_i}(w))
	=\emptyset.
\]	
\end{claim}

\begin{proof}
Assume $v,w\in B(v_i,\beta/2)\cap g^{-m}(B(v_i,\beta/2))$ are $(m,\alpha)$-separated. Arguing by contradiction, suppose that for $v'\eqdef\proj_{D_i}(v)$ it holds $w'\eqdef\proj_{D_i}(w)\in\cC_{i,m}(v')$. As $\cC_{i,m}(v')$ is a proper rectangle, $[v',w']_{R_i}\in\cC_{i,m}(v')$. It follows from \eqref{eq:locpro} that the point
\[\begin{split}
	z
	\eqdef [v,w]
	&= W^\ss_{\varepsilon_{\rm lps}}(v)\cap g^\tau(W^\uu_{\varepsilon_{\rm lps}}(w))\\
	&= W^\ss_{\varepsilon_{\rm lps}}(v)\cap W^\uu_{\varepsilon_{\rm lps}}(w''),
	\quad\text{ where }\,
	\lvert\tau\rvert<\varepsilon_{\rm lps}
	\,\text{ and }\,
	w''\eqdef g^\tau(w),
\end{split}\]
is well defined and $\proj_{D_i}([v,w])=[v',w']_{R_i}$.  As $z\in W^\ss_{\varepsilon_{\rm lps}}(v)$, by Lemma \ref{lem:unifexp} together with Lemma \ref{lem:hamen:uu}, for every $k=0,\ldots,m-1$ it holds
\[
	d(g^k(z),g^k(w''))
	\le cd^\s_{g^k(w'')}(g^k(z),g^k(w''))
	= ce^{-k} d^\s_{w''}(z,w'')
	\le ce^{-k}\varepsilon_{\rm lps}
	< c\varepsilon_{\rm lps}.
\]
On the other hand, as $g^m(v),g^m(w)\in g^m(B(v_i,\beta/2))\cap B(v_i,\beta/2)$, arguing analogously for $g^{-1}$ instead of $g$, it follows that for every $k=0,\ldots,m-1$ it holds
\[
	d(g^k(v),g^k(z))
	< c\varepsilon_{\rm lps}.
\]
Moreover $d(g^k(w''),g^k(w))=d(w'',w)=d(g^\tau(w),w)\le \varepsilon_{\rm lps}$ with $\lvert \tau\rvert\le\varepsilon_{\rm lps}$. It follows that for every $k=0,\ldots,m-1$ it holds
\[
	d(g^k(v),g^k(w))
	\le d(g^k(v),g^k(z))+d(g^k(z),g^k(w''))+d(g^k(w''),g^k(w))
	\le 3c\varepsilon_{\rm lps}
	< \alpha.
\]
But this contradicts the fact that $v,w$ are $\alpha$-separated, proving the claim.
\end{proof}

This finishes the proof of the lemma.
\end{proof}

\subsection{Proof of Proposition \ref{pro:entropybasicset}}\label{ssec:proofprop}

By $G$-invariance, the measure $\mu$ is $g^t$-invariant  for every $t\in\bR$. By Remark \ref{rem:ergodic}, there is $t>0$ such that $\mu$ is $g^t$-ergodic and $g^{-t}$-ergodic. Without loss of generality, we can assume that $t=1$. Let $h\eqdef h_\mu(G)$.

Fix $\varepsilon_{\rm E}\in(0,h/4)$. Fix $r>0$ satisfying
\begin{equation}\label{eq:choicercc}
	r
	<\varepsilon_{\rm E}\cdot\min\Big\{1,\frac{1}{2(h-4\varepsilon_{\rm E})}\Big\}.
\end{equation}
Fix some compact set $Y_0\subset T^1M$ such that $\mu(Y_0)>0$ and $\kappa\in(0,\mu(Y_0)/4)$. 

\subsubsection{Approximate entropy at finite times}

Using notation \eqref{eq:Bowbal} with $f=g^1$, by Brin-Katok's theorem in the noncompact setting (see for example \cite[Theorem 2.7]{Riq:18}), for $\mu$-almost every $x$ it holds
\[
	h
	\le \lim_{\varepsilon\to0}\liminf_{n\to\infty} -\frac1n\log\mu(B_n(x,\varepsilon)). 
\]	
By Egorov's theorem and the fact that $\mu$ is regular, there are a compact set $Y_1\subset Y_0$ satisfying $\mu(Y_1)>\mu(Y_0)-\kappa/4$ and $\varepsilon_{\rm BK}>0$ so that for every $\varepsilon\in(0,\varepsilon_{\rm BK})$ and $x\in Y_1$,
\begin{equation}\label{eq:Egoent}
	h-\frac{\varepsilon_{\rm E}}{3}
	\le \liminf_{n\to\infty} -\frac1n\log\mu(B_n(x,\varepsilon)). 
\end{equation}

\begin{claim}\label{cla:inic}
	For every $\alpha\in(0,\varepsilon_{\rm BK})$ there are $n_{\rm BK}\in\bN$ and a compact set $Y_2\subset Y_1$ satisfying $\mu(Y_2)>\mu(Y_1)-\kappa/4$ such that for every integer $n\ge n_{\rm BK}$ and measurable set $A\subset Y_2$ with $\mu(A)>0$ there exists a $(n,\alpha)$-separated set $E\subset A$ such that
\[
	\frac1n\log\card E 
	\ge h-\frac{2\varepsilon_{\rm E}}{3}-\frac1n\lvert\log\mu(A)\rvert.
\]	
\end{claim}

\begin{proof}
	It follows from \eqref{eq:Egoent} that for every $\alpha \in(0,\varepsilon_{\rm BK})$ there are a compact set $Y_2\subset Y_1$ satisfying $\mu(Y_2)>1-3\kappa/4$ and an integer $n_{\rm BK}\in\bN$ such that for every $n\ge n_{\rm BK}$ and $x\in Y_2$ it holds
\begin{equation}\label{estBK3}
	 \mu(B_n(x,\alpha))
	\le e^{-n(h-2\varepsilon_{\rm E}/3)}. 
\end{equation}

Fix $n\ge n_{\rm BK}$. Let $A\subset Y_2$ be a measurable set with $\mu(A)>0$. We are going to construct an $(n,\alpha)$-separating set $E\subset A$. Choose any point $x_1\in A$. Let $A_1\eqdef A\setminus B_n(x_1,\alpha)$. Inductively, for every $k\ge2$ assuming that $A_{k-1}$ was already constructed and is nonempty, choose $x_k\in A_{k-1}$ and let $A_k\eqdef A_{k-1}\setminus B_n(x_k,\alpha)$.

Let $k\in\bN$ denote the largest index for which $A_k$ is nonempty. By construction, the set $E\eqdef \{x_1,x_2,\ldots,x_k\}$ is $(n,\alpha)$-separated. Now  \eqref{estBK3} implies 
\[
	\mu(A_k)
	= \mu(A_{k-1})-\mu(B_n(x_k,\alpha))
	= \mu(A)-\sum_{i=1}^k\mu(B_n(x_k,\alpha))
	\ge \mu(A)-ke^{-n(h-2\varepsilon_{\rm E}/3)}.
\]
Hence, it holds $k\ge\lceil \mu(A)e^{n(h-2\varepsilon_{\rm E}/3)}\rceil$ and 
\[
	\frac1n\log\card E
	= \frac1n\log k
	\ge h -\frac{2\varepsilon_{\rm E}}{3}-\frac1n\lvert\log\mu(A)\rvert,
\]
proving the claim.
\end{proof}

Once $\alpha\in(0,\varepsilon_{\rm BK})$ is fixed, in any measurable set of positive measure, for sufficiently large $n$ we can find an $(n,\alpha)$-separated set of points whose cardinality is close to $e^{nh}$. We now choose the maximal size of cross sections and rectangles and some finite partition to construct a horseshoe. Note that the exponential growth of the cardinality of separate points in Claim \ref{cla:inic} occurs in every partition element (in fact in every measurable set of positive measure). The principle idea in \cite[Supplement S.5]{KatHas:95} is to pick \emph{one} element with the largest cardinality, which still will be of order $e^{nh}$ for $n$ large. Then, we consider a corresponding local section which contains this partition element and build the horseshoe as the invariant set of a certain return-map to that section.

\subsubsection{Fixing quantifiers and local cross sections}

Fix $\alpha\in(0,\varepsilon_{\rm BK})$, let $n_{\rm BK}$ and $Y_2$ as provided by Claim \ref{cla:inic}.  Without loss of generality, we can assume that $n_{\rm BK}$ is so large  that for every $n\ge n_{\rm BK}$ it holds
\[
	e^{-n\varepsilon_{\rm E}/3}
	\le 1-\frac{3\kappa}{4}
	\quad\text{ and }\quad
	n
	<e^{n\varepsilon_{\rm E}}.
\]
Apply Lemma \ref{lem:rec1} to $Y_2$ and $\alpha$. For $\varepsilon_{\rm sec}\in(0,\alpha/(3c))$ small enough this lemma provides numbers $\beta_{\rm sec}>0$ and $M\in\bN$ and a finite collection of sections $\{D_i\}_{i=1}^\ell$ of size smaller than $\varepsilon_i\in(0,\varepsilon_{\rm sec})$. Moreover, fixing some
\[
	\beta\in(0,\beta_{\rm sec}).
\]
there are points $v_i\in Y_2$ such that
\[
	Y_2
	\subset\bigcup_iB(v_i,\frac\beta2)
	\quad\text{ and }\quad
	B(v_i,\frac\beta2)
	\subset B(v_i,\beta)
	\subset \widehat R_i\eqdef g^{(-\varepsilon_i,\varepsilon_i)}(R_i), 
\] 
where $R_i$ is a proper rectangle in $D_i$.

\subsubsection{Fixing a partition}
Fix a finite partition $\cP=\{P_1,\ldots,P_\ell\}$ of $Y_0$ of diameter at most $\beta/4$.  Denoting by $\cP(x)$ the partition element which contains $x$. Hence, with this choice,
\[
	\cP(v_i)
	\subset B(v_i,\frac\beta2)
	\subset B(v_i,\beta)
	\subset \widehat R_i. 
\]

\subsubsection{Choose almost-uniformly returning points}
By Birkhoff ergodic theorem (relative to the $g^1$-ergodic measure $\mu$), for $\mu$-almost every $x\in Y_2$ for every $i$ it holds
\[
	\lim_{n\to\infty}\frac1n\card\big\{k\in\{0,\ldots,n-1\}\colon g^k(x)\in Y_2\cap\cP(v_i)\big\}
	=\mu(Y_2\cap\cP(v_i)).
\]
Let $\varepsilon_{\rm B}\in(0,r)$. By Egorov's theorem, there is an integer $n_{\rm B}\in\bN$ and a compact set $Y_3\subset Y_2$ 
such that
\begin{equation}\label{eq:choiceY3}
	\mu\big(Y_3\big)
	> \mu(Y_2)-\frac\kappa 4
	> 1-\kappa.
\end{equation}
 such that for every index $i$,  point $x\in Y_3\cap\cP(v_i)$, and $n\ge n_{\rm B}$ it holds
\begin{equation}\label{eq:choicen0}\begin{split}	
	&\big\lvert \card\{k\in\{0,\ldots,n-1\}\colon g^k(x)\in Y_2\cap \cP(v_i)\} 
		-n\mu(Y_2\cap \cP(v_i))\big\rvert
	\le n\varepsilon_{\rm B}.
\end{split}\end{equation}
Without loss of generality, we can assume that $n_{\rm B}$ is so large that 
\begin{equation}\label{eq:choicen1}
	n_{\rm B}r\big(\mu(Y_2)-3\varepsilon_{\rm B}\big)
	\ge1.
\end{equation}
Hence, for every $i$, $x\in Y_3\cap \cP(v_i)$, and $n\ge n_{\rm B}$,  it holds
\[\begin{split}
	\card&\big\{k\in\{n,\ldots,n(1+r)-1\}\colon g^k(x)\in Y_2\cap \cP(v_i)\big\}\\
\text{\small by \eqref{eq:choicen0}}\quad
	&\ge n(1+r)\big(\mu(Y_2\cap \cP(v_i))-\varepsilon_{\rm B}\big)-(n-1)\big(\mu(Y_2\cap \cP(v_i))+\varepsilon_{\rm B}\big)\\
	&= (nr+1)\mu(Y_2\cap \cP(v_i))-(2n+nr-1)\varepsilon_{\rm B}\\
\text{\small using $\varepsilon_{\rm B}<r$ }\quad
	&> nr\big(\mu(Y_2\cap \cP(v_i))-3\varepsilon_{\rm B}\big)\\
\text{\small by \eqref{eq:choicen1} and $n\ge n_{\rm B}$}\quad
	&\ge 1.
\end{split}\]
In other words, there exists $k=k(x)\in\bN$ satisfying $k\in\{n,\ldots,n(1+r)-1\}$ such that
\[
	g^k(x)\in \cP(x)= \cP(v_i),
\]
that is, the point returns to its partition element. 

Apply now the Claim \ref{cla:inic} to the set $A=Y_3$. Fix $n\in\bN$ satisfying
\begin{equation}\label{eq:choiceY3b}
	n\ge \max\Big\{n_{\rm BK},M,n_{\rm B},
		\frac{3}{\varepsilon_{\rm E}}\lvert\log(1-\kappa)\rvert,
		\frac{1}{\varepsilon_{\rm E}}\log\ell,
		\frac{h}{2\varepsilon_{\rm E}}+\frac32,
		\frac{\varepsilon_{\rm sec}}{r}
		\Big\}.
\end{equation}
By Claim \ref{cla:inic}, there is a $(n,\alpha)$-separated set $E\subset Y_3$ such that
\begin{equation}\label{eq:cardEest}\begin{split}
	\frac{1}{n}\log\card E
	&\ge h-\frac{2\varepsilon_{\rm E}}{3} - \frac1n\lvert\log\mu(Y_3)\rvert\\
	\text{\small by \eqref{eq:choiceY3}}\quad	
	&\ge h-\frac{2\varepsilon_{\rm E}}{3} - \frac1n\lvert\log(1-\kappa)\rvert	\\
	\text{\small by \eqref{eq:choiceY3b}}\quad
	&> h-\varepsilon_{\rm E}.
\end{split}\end{equation}

\subsubsection{Picking points with the same return time to the same partition element}
For each $k\in\{n,\ldots,n(1+r)-1\}$ let
\[
	E_k
	\eqdef \{x\in E\colon g^k(x)\in\cP(x)\}
\]
be the set of points in $E$ that return to their partition element at the same time $k$.
Let 
\begin{equation}\label{eq:choicem}
	m\in\{n,\ldots,n(1+r)-1\}
\end{equation}	
 be an index $k$ for which $\card E_k$ is maximal. Since
\[
	\card E
	= \sum_{k=n}^{n(1+r)-1}\card E_k,
\]
it follows $\card E\le rn\card E_m$. Together with $rn<e^{rn}$ and \eqref{eq:cardEest}, it follows
\[
	\card E_m
	\ge \frac{1}{rn}\card E
	> e^{-rn}\card E
	> e^{n(h-\varepsilon_{\rm E})-rn}.
\]
Pick from the partition $\cP$ an element $P_i$ for which $\card (E_m\cap P_i)$ is maximal. As $\cP$ has $\ell$ elements, 
\begin{equation}\label{eq:estcardN}\begin{split}
	\card(E_m\cap P_i)
	&\ge \frac{1}{\ell}\card E_m
	\ge \frac1\ell e^{n(h-\varepsilon_{\rm E})-rn}\\
	\text{\small by \eqref{eq:choiceY3b}}\quad
	&> e^{n(h-2\varepsilon_{\rm E}-r)}\\
	\text{\small by \eqref{eq:choicercc}}\quad
	&> e^{n(h-3\varepsilon_{\rm E})}.
\end{split}\end{equation}
For what is below, fix this rectangle $R_i$ in the cross section $D_i$.  Fix some enumeration 
\[
	F=\{x_1,\ldots,x_N\}
	\eqdef E_m\cap P_i.
\]	

\subsubsection{Building the basic set}
 Note that all points in $F$ are $(n,\alpha)$-separated and hence $(m,\alpha)$-separated. Recall that  it holds
\[
	F\cap g^m(F)
	\subset \cP(v_i)
	\subset B(v_i,\frac\beta2)
	\subset B(v_i,\beta)
	\subset \widehat R_i
	= g^{(-\varepsilon_i,\varepsilon_i)}(R_i).	
\] 
We now are ready to build the basic set around the set $F$ in the neighborhood of $v_i$ and the points $x_1,\ldots,x_N$. 

Consider the continuous and injective return map $T_{i,m}\colon\widehat D_i^m\to\widehat D_i^m$
 as provided by Lemma \ref{lem:rec1}. 
Denote the connected component containing $x_\ell\in F$, $\ell=1,\ldots,N$, by
\[
	U_\ell
	\eqdef \cC_{i,m}(\proj_{D_i}(x_\ell)).
\]
By Lemma \ref{lem:rec1}, this are proper rectangles in $D_i$ which are pairwise disjoint. Let
\[
	S_\ell
	\eqdef T_{i,m}(U_\ell)
	\subset R_i,
\]
which, by the above, are also pairwise disjoint and closed. Let 
\[
	\Gamma
	\eqdef \bigcap_{k\in\bZ}T_{i,m}^k(U_1\cup\ldots\cup U_N),
\]
and note that this is a compact set. As $U_\ell$ are pairwise disjoint, for every $x\in\Gamma$, there is a unique sequence $\underline i=(\ldots i_{-1}i_0i_1\ldots)\eqdef\pi(x)\in\Sigma_N$ of indices such that
\[
	T_{i,m}^k(x)
	\in U_{i_k}
\]
for every $k\in\bZ$. This defines $\pi\colon\Gamma\to\Sigma_N$. By definition, $\pi\circ T_{i,m}=\sigma\circ \pi$. By classical arguments, $\Gamma$ is a topological horseshoe and $\pi$ is injective and onto.

\begin{claim}\label{cla:picont}
	The map $\pi^{-1}\colon\Sigma_N\to\Gamma$ is continuous. 
\end{claim}

\begin{proof}
Let $\varepsilon>0$ and take $k\in\bN$ such that $Ce^{-k(m-\varepsilon_{\rm sec})}<\varepsilon$, where $C=C(\varepsilon_{\rm sec})$ is the constant of Lemma \ref{lem:closegeodesics}. Let $\delta=1/2^k$. Then, whenever $\underline i$, $\underline i'$ verify $d_\Sigma(\underline i,\underline i')<\delta$, that is $i_j=i'_j$ for all $-k\leq j\leq k$, we get that for $x=\pi^{-1}(\underline i)$ and $x'=\pi^{-1}(\underline i')$, the orbits $g^t(x)$ and $g^t(x')$ are $\varepsilon_{\rm sec}$-close for time at least $2k(m-\varepsilon_{\rm sec})$. Using Lemma \ref{lem:closegeodesics} on the quotient manifold, it follows that $d(x,y)\leq Ce^{-k(m-\varepsilon_{\rm sec})}<\varepsilon$, which proves the claim. 
\end{proof}

Let 
\[
	B
	\eqdef \bigcup_{\lvert t\rvert\le\max T_{i,m}}g^t(\Gamma).
\]
It follows from Claim \ref{cla:picont} that $B$ is a $G$-invariant compact set of topological dimension 1.

Let us now estimate the topological entropy of $G$ on $B$.
Note that $G|_B$ can be seen as the suspension flow of the map $T_{i,m}|_\Gamma$ under the return time defined as in \eqref{eq:estTim}. Clearly,
\[
	h_{\rm top}(T_{i,m},\Gamma)
	= \log N.
\]
Consider on $\Gamma$ the push-forward of the $(\frac1N,\ldots,\frac1N)$-Bernoulli measure by $\pi$ and denote by $\mu$ its suspension. Recall that $\mu$ is a $G$-invariant Borel probability measure. As by \eqref{eq:estTim} the return time is $\varepsilon_i$-close to $m$, by Abramov's formula, it holds
\[
	h_{\rm top}(G,\Gamma)
	\ge h_\mu(G)
	\ge \frac{\log N}{m+\varepsilon_i},
\]
 Together with \eqref{eq:estcardN} and \eqref{eq:choicem}, it follows
\[\begin{split}
	\frac{\log N}{m+\varepsilon_i}
	&\ge \frac{n(h-3\varepsilon_{\rm E})}{n(1+r)+\varepsilon_i}\\
	\text{\small recalling that $\varepsilon_i<\varepsilon_{\rm sec}$ and using \eqref{eq:choiceY3b}}
	&\ge \frac{h-3\varepsilon_{\rm E}}{1+2r}\\
	\text{\small by \eqref{eq:choicercc}}\quad
	&\ge h-4\varepsilon_{\rm E}.
\end{split}\]
This implies the assertion of the proposition. \qed

\section{Exceptional sets for suspension flows}\label{sec:susflo}

In this section, we prove Theorem \ref{thepro:susflo}. Let us first recall the definition of a suspension flow. Given $N\in\bN$, consider the shift space $\Sigma=\{1,\ldots,N\}^\bZ$. Given $\alpha>0$, consider a continuous \emph{height function} $\tau \colon \Sigma\to(0,\alpha]$ and define the space
\[
	\Sigma(\sigma,\tau)
	\eqdef \big(\Sigma\times[0,\alpha]\big)_{\sim},
\]
as the quotient space of $\Sigma\times\bR_{\ge0}$ modulo the equivalence relation $\sim$ that identifies $(\underline i,s)$ with $(\sigma(\underline i),s-\tau(\underline i))$ for all $\underline i\in\Sigma$ and $s\ge\tau(\underline i)$. This space is compact and metrizable (see \cite[Section 2]{BowWal:72}). 
The \emph{suspension flow} $F=(f^t)_{t\in\bR}$ of the left-shift $\sigma\colon\Sigma\to\Sigma$ \emph{under $\tau$} is the map $F\colon\Sigma(\sigma,\tau)\times\bR\to\Sigma(\sigma,\tau)$ defined by 
\[
	f^t(\underline i,s)
	\eqdef \begin{cases}
	(\underline i,s+t)&\text{if }0\le s+t<\tau(\underline i),\\
	(\sigma(\underline i),s+t-\tau(\underline i))&\text{if }s+t=\tau(\underline i).
	\end{cases}
\]

\begin{proof}[\bf{Proof of Theorem \ref{thepro:susflo}}]
Note that the suspension flow $F$ is expansive and hence the map $\mu\mapsto h_\mu(F)$ is upper semi-continuous. Hence, there is some ergodic measure of maximal entropy $\mu$, $h_\mu(F)=h_{\rm top}(F)$.  Considering the natural projection $\pi_1\colon\Sigma(\sigma,\tau)\to\Sigma$ to the first coordinate and let $\nu_{\rm max}\eqdef(\pi_1)_\ast\mu\in\cM_{\rm erg}(\sigma)$.

Recall that for every $\nu\in\cM(\sigma)$
\[
	\mu_\nu
	\eqdef \frac{1}{(\nu\times m)(\Sigma(\sigma,\tau))}(\nu\times m)|_{\Sigma(\sigma,\tau)}
\]
defines a $F$-invariant Borel probability measure. Moreover, $\nu\mapsto\mu_\nu$ is a bijection between $\cM(\sigma)$ and $\cM(F)$, and it holds
\begin{equation}\label{eq:entropysusflo}
	h_{\mu_\nu}(F)
	=\frac{h_\nu(\sigma)}{\int \tau\,d\nu}.
\end{equation}
Hence, together with our hypothesis on $A$, it holds
\begin{equation}\label{eq:using}
	h^\ast_{\rm top}(f^1,A)
	< h_{\rm top}(f^1)
	= h_\mu(F)
	= \frac{h_{\nu_{\rm max}}(\sigma)}{R},
	\quad\text{where}\quad
	R\eqdef \int\tau\,d\nu_{\rm max}.
\end{equation}
Choose $\varepsilon>0$ sufficiently small such that 
\begin{equation}\label{eq:suchthat}
	h^\ast_{\rm top}(f^1,A)
	< \frac{R-\varepsilon}{R+\varepsilon} \frac{h_{\nu_{\rm max}}(\sigma)-\varepsilon}{R}.
\end{equation}
Apply Proposition \ref{pro:SFT} to $\nu_{\rm max}$ and $\varepsilon$ and consider the corresponding SFT $\Xi\subset\Sigma$ which hence satisfies
\begin{equation}\label{eq:using2}
	h_{\nu_{\rm max}}(\sigma)-\varepsilon
	\le h_{\rm top}(\sigma|_\Xi)
	\quad\text{ and }\quad
	\left\lvert\int\tau\,d\nu_{\rm max} -\int\tau\,d\lambda\right\rvert
	 \le\varepsilon
\end{equation}
In the following, we consider the corresponding suspension flow on the suspension space $\Xi(\sigma,\tau)=(\Xi\times[0,\alpha])_\sim$ which can be considered as a subset of $\Sigma(\sigma,\tau)$. Let $W\eqdef A\cap\Xi(\sigma,\tau)$ and observe that by Lemma \ref{lem:propTho} (ii) and \eqref{eq:suchthat}, it holds
\begin{equation}\label{eq:relatt}
	h^\ast_{\rm top}(f^1,W)
	\le h^\ast_{\rm top}(f^1,A)
	< \frac{R-\varepsilon}{R+\varepsilon} \frac{h_{\nu_{\rm max}}(\sigma)-\varepsilon}{R}.
\end{equation}

For the following claim, recall notation \eqref{notationMZ}.

\begin{claim}\label{cla:sameproj}
	It holds 
\[
	\sup_{\nu\in\cM^{\pi_1(W)}(\sigma)}\frac{h_\nu(\sigma)}{\int\tau\,d\nu}
	= 
	\sup_{\mu\in\cM^W(f^1)}h_\mu(f^1).
\]	
\end{claim}

\begin{proof}
The natural projection $\pi_1$ induces a continuous push forward $(\pi_1)_\ast\colon\cM(f^1)\to\cM(\sigma)$. Hence if $\mu'\in\cV(f^1,X)$ then $(\pi_1)_\ast\mu'\in\cV(\sigma,\pi_1(X))$, which implies 
\[
	\cM^{\pi_1(W)}(\sigma)
	\supset(\pi_1)_\ast(\cM^W(f^1)).
\]	 
Together with \eqref{eq:entropysusflo} this implies the inequality $\ge$ in the claim. 

On the other hand, let $\xi\in \pi_1(W)$ and $s\in[0,\alpha)$ such that $X=(\xi,s)\in W$ and let $\nu\in\cM^{\pi_1(W)}(\sigma)$ be the limit point of the sequence of probability measures $(\delta_{\xi,n})_n$ as defined in \eqref{eq:deltadef}   (with respect to the shift map $\sigma$). To simplify notation, assume that this sequence converges as $n\to\infty$ (otherwise pass to some subsequence). Consider a subsequence $(k_n)_n$ so that $\pi_1(f^{k_n}(X))=\sigma^n(\xi)$ and let $\mu'$ be some limit point of $(\delta_{X,k_n})_n$ (relative to the space of probability measures in $\Sigma(\sigma,\tau)$, possibly for some subsequence). Then $\mu'$ is $f^1$-invariant and $(\pi_1)_\ast\mu'=\nu$. The latter implies $h_{\mu'}(f^1)=h_\nu(\sigma)/\int\tau\,d\nu$. This implies the inequality $\le$  in the claim.
\end{proof}

Denote $Z\eqdef \pi_1(W)$. Observe that $Z\subset\Xi$. Note that $\xi\in Z$ implies that every limit measure $\lambda$ of $\delta_{\xi,n}$ (relative to the shift map $\sigma\colon\Sigma\to\Sigma$ and the weak$\ast$ topology in the space of probability measures in $\Sigma$) is in $\cM(\sigma|_\Xi)$. Hence, by our choice of SFT $\Xi$ and \eqref{eq:using}, 
\begin{equation}\label{eq:variationtau}
	R-\varepsilon
	\le \int\tau\,d\lambda \le R+\varepsilon.
\end{equation}
By Claim \ref{cla:sameproj} and definition \eqref{def:entstar} it follows
\[\begin{split}
	\sup_{\nu\in\cM^Z(\sigma)}\frac{h_\nu(\sigma)}{\int\tau\,d\nu}
	&= \sup_{\mu\in\cM^W(f^1)}h_\mu(f^1)
	= h^\ast_{\rm top}(f^1,W)\\
		\text{\small{with \eqref{eq:relatt}}}\quad
	&< \frac{R-\varepsilon}{R+\varepsilon} \frac{h_{\nu_{\rm max}}(\sigma)-\varepsilon}{R}
\end{split}\]
This together with \eqref{eq:variationtau}, this implies
\[\begin{split}
	\sup_{\nu\in\cM^Z(\sigma)}{h_\nu(\sigma)}
	&\le (R+\varepsilon)\cdot\frac{R-\varepsilon}{R+\varepsilon} \frac{h_{\nu_{\rm max}}(\sigma)-\varepsilon}{R} \\
	&< h_{\nu_{\rm max}}(\sigma)-\varepsilon\\
		\text{\small by \eqref{eq:using2}}\quad
	&\le h_{\rm top}(\sigma|_\Xi).
\end{split}\]
Together with  Lemma \ref{lem:propTho2} (iii), it follows
\begin{equation}\label{eq:gelo}
	h_{\rm top}(\sigma|_\Xi,Z)
	\le h^\ast_{\rm top}(\sigma|_\Xi,Z)
	= \sup_{\nu\in\cM^Z(\sigma)}{h_\nu(\sigma)}
	<h_{\rm top}(\sigma|_\Xi).
\end{equation}
Hence, by \cite[Section 3 Theorem 1]{Dol:97}, this implies
\[
	h_{\rm top}(\sigma|_\Xi, I_{\sigma|\Xi}(Z))
	= h_{\rm top}(\sigma|_\Xi)
\]
By Lemma \ref{lem:propTho2} (iii), it holds
\[
	h_{\rm top}(\sigma|_\Xi, I_{\sigma|\Xi}(Z))
	\le h^\ast_{\rm top}(\sigma|_\Xi, I_{\sigma|\Xi}(Z))
\]
and therefore
\begin{equation}\label{eq:conDol}
       h_{\rm top}(\sigma|_\Xi, I_{\sigma|\Xi}(Z))
       = h^\ast_{\rm top}(\sigma|_\Xi, I_{\sigma|\Xi}(Z))
       = h_{\rm top}(\sigma|_\Xi).
\end{equation}

\begin{claim}\label{cla:sameproj2}
	It holds $I_{\sigma|\Xi}(\pi_1(W))\subset \pi_1(I_{f^1}(W))$.
\end{claim}

\begin{proof}
Let $X\not\in I_{f^1}(W)$ with $f^{k_n}(X)\to Y\in W$ for some subsequence $k_n\to\infty$. Let $\xi=\pi_1(X)$. By continuity of $\pi_1$, it follows that there is some subsequence $\ell_m\to\infty$ so that $\sigma^{\ell_m}(\xi)=(\pi_1\circ f^{k_n})(X)\to\pi_1(Y)\in\pi_1(W)=Z$ and hence $\xi\not\in I_{\sigma|\Xi}(Z)$.
\end{proof}

Let $Z\eqdef \pi_1(W)$. Claim \ref{cla:sameproj2} implies
\[
	\cM^{I_{\sigma|\Xi}(\pi_1(W))}(\sigma)
	=\cM^{I_{\sigma|\Xi}(Z)}(\sigma)
	\subset \cM^{\pi_1(I_{f^1}(W))}(\sigma).
\]
This together with Claim \ref{cla:sameproj} implies
\[
	\sup\left\{\frac{h_\lambda(\sigma)}{\int\tau\,d\lambda}\colon 
		\lambda\in\cM^{I_{\sigma|\Xi}(Z)}(\sigma)\right\}
	\le\sup\Big\{h_\mu(f^1)\colon \mu\in \cM^{I_{f^1}(W)}(f^1)\Big\}.
\]
Using \eqref{eq:conDol} together with \eqref{eq:variationtau}, this implies
\begin{equation}\label{eq:auto}
	 \frac{h_{\rm top}(\sigma|_\Xi)}{R+\varepsilon}
	 = \frac{h^\ast_{\rm top}\big(\sigma|_\Xi,I_{\sigma|\Xi}(Z)\big)}{R+\varepsilon}
	\le h^\ast_{\rm top}(f^1,I_{f^1}(W)).
\end{equation}
By  \eqref{eq:using} and \eqref{eq:using2}, it holds
\[
	h_\mu(f^1)
	 \frac{R}{h_{\nu_{\rm max}}(\sigma)}
	\frac{h_{\nu_{\rm max}}(\sigma)-\varepsilon}{R+\varepsilon}
	= 1\cdot\frac{h_{\nu_{\rm max}}(\sigma)-\varepsilon}{R+\varepsilon}
	\le \frac{h_{\rm top}(\sigma|_\Xi)}{R+\varepsilon}.
\]
Recalling again that $\mu$ was a measure of maximal entropy $h_\mu(f^1)=h_{\rm top}(F)$, together with \eqref{eq:auto}, Lemma \ref{lem:propTho} (ii) applied to $I_{f^1}(W)\subset\Sigma(\sigma,\tau)$, and Lemma \ref{lem:propTho2} (iv), we get
\[
	h_{\rm top}(F)
	\frac{R(h_{\nu_{\rm max}}(\sigma)-\varepsilon)}{h_{\nu_{\rm max}}(\sigma)
	(R+\varepsilon)}
	\le h^\ast_{\rm top}(f^1,I_{f^1}(W))
	\le h^\ast_{\rm top}(f^1,\Sigma(\sigma,\tau))
	 = h_{\rm top}(f^1,\Sigma(\sigma,\tau)).
\]
As $\varepsilon>0$ was arbitrary, this proves the theorem.
\end{proof}

We apply now Theorem \ref{thepro:susflo} to prove the following result.

\begin{proposition}\label{pro:basiexce}
	Let $B\subset T^1M$ be a one-dimensional basic set for the geodesic flow $G$.
	For every Borel set $A\subset T^1M$ satisfying $h^\ast_{\rm top}(G,A\cap B)<h_{\rm top}(G|_B)$ it holds
\[
	h^\ast_{\rm top}(G|_B,I_{G|B}(A))
	= h_{\rm top}(G|_B).
\]		
\end{proposition}

\begin{proof}
Note that the flow $G$ restricted to $B$ is topologically conjugate to a suspension flow $F$ of the shift map $\sigma$ on some subshift of finite type $\Sigma\subset\{1,\ldots,N\}^\bZ$ under some continuous function $\tau\colon\Sigma\to[0,\alpha)$, $\alpha>0$ \cite[Theorem 10]{BowWal:72}. Denote by $p\colon \Sigma(\sigma,\tau)\to B$ the corresponding conjugating homeomorphism so that for every $t\in\bR$,
\[
	p\circ f^t=g^t\circ p.
\]
Consider the Borel sets 
\[
	W
	\eqdef p^{-1}(A\cap B),
	\quad
	Z
	\eqdef (\pi_1\circ p^{-1})(A\cap B)
	= \pi_1(W).
\]	
By Lemma \ref{lem:propTho2} (ii), it follows 
\begin{equation}\label{eq:ent1}
	h_{\rm top}^\ast({g^1}|_B,A\cap B)
	= h_{\rm top}^\ast({f^1},p^{-1}(A\cap B))
	= h_{\rm top}^\ast({f^1},W).
\end{equation}
By conjugation, it holds
\[
	h_{\rm top}(g^1|_B)
	= h_{\rm top}(f^1).
\]
Hence, it follows from our hypothesis that
\[
	h_{\rm top}^\ast({f^1},W)
	<h_{\rm top}(f^1).
\]
By Theorem \ref{thepro:susflo}, it holds
\[
	h_{\rm top}^\ast(f^1,I_{f^1}(W))
	= h_{\rm top}(f^1).
\]
It remains to check that
\[
	p(I_{f^1}(W))
	= I_{g^1|_B}(p(W))
	= I_{g^1|_B}(A\cap B)
	= I_{g^1|_B}(A).
\]
Applying again Lemma \ref{lem:propTho2} (ii), we get $h^\ast_{\rm top}(G|_B,I_{g^1|B}(A))= h_{\rm top}(G|_B)$. Note that $I_{g^1|B}(A)= \bigcup_{t\in [0,1)} g^t I_{G|B}(A)$ so a measure $\mu\in \mathcal{M}^{I_{g^1|B}(A)}(g^1)$ defines a measure $(g^{-t})_\ast\mu\in\mathcal{M}^{I_{G|B}(A)}(g^1)$. Since $\mu$ can be chosen $G$-invariant, the variational principle implies 
\[
      h^\ast_{\rm top}(G|_B,I_{g^1|B}(A)) 
      = h^\ast_{\rm top}(G|_B,I_{G|B}(A)), 
\]
and the claim follows.
\end{proof}

\section{Proofs of Theorem \ref{the:1} and Corollaries \ref{cor:1} and \ref{cor:2}}\label{sec:proofs}

\begin{proof}[Proof of Theorem \ref{the:1}]
Consider a Borel set $A\subset T^1M$ satisfying
\[
	h^\ast_{\rm top}(G,A)
	< h_{\rm top}(G).
\]
By \eqref{eq:timetmap} and the variational principle \eqref{eq:vpflow} together with the definition \eqref{def:entstar}, 
\[
	h^\ast_{\rm top}(G,A)
	= \sup_{\mu\in\cM^A(g^1)}h_\mu(g^1)
	< \sup_{\mu\in\cM(g^1)}h_\mu(g^1)
	= h_{\rm top}(G).
\]
Let $\mu\in\cM(g^1)$ satisfying
\[
	h^\ast_{\rm top}(G,A)
	< h_\mu(g^1).
\]
Without loss of generality (otherwise, choose some ergodic component having this property), we can assume that $\mu$ is ergodic. Fix $\varepsilon>0$ such that
\[
	h^\ast_{\rm top}(G,A)
	< h_\mu(g^1)-2\varepsilon.
\]
By Proposition \ref{pro:entropybasicset}, there is a basic set $B\subset T^1M$ of topological dimension 1 satisfying
\[
	h_\mu(g^1)-\varepsilon
	< h_{\rm top}(G|_B).
\] 
By Lemma \ref{lem:propTho} (ii), it holds 
\[
	h^\ast_{\rm top}(G|_B,A\cap B)
	\le h^\ast_{\rm top}(g^1,A)
	< h_{\rm top}(G|_B).
\]
Proposition \ref{pro:basiexce} now implies 
\[
	h_{\rm top}^\ast(G|_B,I_{G|B}(A))
	= h_{\rm top}(G|_B).
\]

Since $B$ is $G$-invariant, we have $I_{G|B}(A)\subset I_{G}(A)$, therefore
\[ 
             h_\mu(g^1)-\varepsilon<h_{\rm top}^\ast(G,I_{G}(A))
\]
As $\varepsilon>0$ was arbitrary, from the variational principle \eqref{eq:vpflow} we conclude that $h_{\rm top}(G)=h_{\rm top}^\ast(G,I_{G}(A))$.
\end{proof}

\begin{proof}[Proof of Corollary \ref{cor:1}]
By assumption the derivatives of the sectional curvatures are uniformly bounded, so there is at most one measure maximizing entropy for $G$ (see \cite[Theorem 1]{OP:04}). If $h^{\ast}_{\rm top}(G,K)=h_{\rm top}(G)$, then the measure of maximal entropy of $(K,G|_{K})$ is the unique measure of maximal entropy for $G$. This measure is known to be fully supported on $\Omega$, a contradiction since $\Omega$ is non-compact. Therefore $h^{\ast}_{\rm top}(G,K)<h_{\rm top}(G)$, and by Theorem \ref{the:1} we obtain the desired equality.
\end{proof}

\begin{proof}[Proof of Corollary \ref{cor:2}]
Let us first show that if $\Gamma'$ is divergent and $\Omega_N \neq \Omega$, then $\delta_{\Gamma'} < \delta_{\Gamma}$. Indeed, the condition $\Omega_N \neq \Omega$ is equivalent to have $L(\Gamma')\neq L(\Gamma)$, where $L(H)$ denotes the limit set of any discrete subgroup of isometries $H$ of $\widetilde{M}$. Then \cite[Proposition 2]{DOP:00} implies $\delta_{\Gamma'} < \delta_{\Gamma}$. 
Note that by \cite[Theorem 1]{OP:04}, it holds $h_{\rm top}(G)=\delta_\Gamma$ and $h_{\rm top}(G,N)=\delta_{\Gamma'}$. 
Hence, it holds $h_{\rm top}(G,N)<h_{\rm top}(G)$. This together with Theorem \ref{the:1} ends the proof.
\end{proof}

\begin{proof}[Proof of Scholium \ref{sch:1}]
	By the variational principle \eqref{eq:vpflow}, for every $\varepsilon>0$ there exists an ergodic measure $\mu\in\cM_{\rm erg}(G)$ satisfying $h_\mu(G)\ge h_{\rm top}(G,T^1M)-\varepsilon$. By Proposition \ref{pro:entropybasicset}, there exists a basic (and hence, in particular bounded, set $B\subset T^1M$ satisfying $h_{\rm top}(G,B)\ge h_\mu(G)-\varepsilon$. As $\varepsilon$ was arbitrary, the supremum taken over all basic sets equals $h_{\rm top}(G,T^1M)$. The remaining equalities in the scholium follow from Lemma \ref{lem:propTho} item (ii) and Lemma \ref{lem:propTho2} item (iii). 
\end{proof}

\subsection{Final comments} Most of the results in this work are not proved using the geometric structure of the dynamical system, but its hyperbolic structure. Nevertheless, the estimates from Section \ref{sec:geodynneg} are crucially used, for instance in the construction of suitable local cross sections (Lemma \ref{lem:rec1}). It would be interesting to know how general the results in this work can be for more general classes of hyperbolic flows in noncompact manifolds.

\bibliographystyle{alpha}
\bibliography{bib}

\begin{thebibliography}{DOP00}

\bibitem[Abr59]{Abr:59}
L.~M. Abramov.
\newblock On the entropy of a flow.
\newblock {\em Dokl. Akad. Nauk SSSR}, 128:873--875, 1959.

\bibitem[Bea95]{Bea:83}
Alan~F. Beardon.
\newblock {\em The geometry of discrete groups}, volume~91 of {\em Graduate
  Texts in Mathematics}.
\newblock Springer-Verlag, New York, 1995.
\newblock Corrected reprint of the 1983 original.

\bibitem[BK83]{BriKat:83}
M.~Brin and A.~Katok.
\newblock On local entropy.
\newblock In {\em Geometric dynamics ({R}io de {J}aneiro, 1981)}, volume 1007
  of {\em Lecture Notes in Math.}, pages 30--38. Springer, Berlin, 1983.

\bibitem[Bow73]{Bow:73}
Rufus Bowen.
\newblock Topological entropy for noncompact sets.
\newblock {\em Trans. Amer. Math. Soc.}, 184:125--136, 1973.

\bibitem[Bow08]{Bow:08}
Rufus Bowen.
\newblock {\em Equilibrium states and the ergodic theory of {A}nosov
  diffeomorphisms}, volume 470 of {\em Lecture Notes in Mathematics}.
\newblock Springer-Verlag, Berlin, revised edition, 2008.
\newblock With a preface by David Ruelle, Edited by Jean-Ren\'{e} Chazottes.

\bibitem[BW72]{BowWal:72}
Rufus Bowen and Peter Walters.
\newblock Expansive one-parameter flows.
\newblock {\em J. Differential Equations}, 12:180--193, 1972.

\bibitem[CG16]{CaGe:16}
Sara Campos and Katrin Gelfert.
\newblock Exceptional sets for nonuniformly expanding maps.
\newblock {\em Nonlinearity}, 29(4):1238--1256, 2016.

\bibitem[CG19]{CaGe:19}
Sara Campos and Katrin Gelfert.
\newblock Exceptional sets for nonuniformly hyperbolic diffeomorphisms.
\newblock {\em J. Dynam. Differential Equations}, 31(2):979--1004, 2019.

\bibitem[CLT20]{ConLafTho:20}
David Constantine, Jean-Fran\c{c}ois Lafont, and Daniel~J. Thompson.
\newblock Strong symbolic dynamics for geodesic flows on {${\rm CAT}(-1)$}
  spaces and other metric {A}nosov flows.
\newblock {\em J. \'{E}c. polytech. Math.}, 7:201--231, 2020.

\bibitem[Dan86]{Dan:86}
S.~G. Dani.
\newblock Bounded orbits of flows on homogeneous spaces.
\newblock {\em Comment. Math. Helv.}, 61(4):636--660, 1986.

\bibitem[Din70]{Din:70}
E.~I. Dinaburg.
\newblock A correlation between topological entropy and metric entropy.
\newblock {\em Dokl. Akad. Nauk SSSR}, 190:19--22, 1970.

\bibitem[Dol97]{Dol:97}
D.~Dolgopyat.
\newblock Bounded orbits of {A}nosov flows.
\newblock {\em Duke Math. J.}, 87(1):87--114, 1997.

\bibitem[DOP00]{DOP:00}
Fran\c{c}oise Dal'bo, Jean-Pierre Otal, and Marc Peign\'{e}.
\newblock S\'{e}ries de {P}oincar\'{e} des groupes g\'{e}om\'{e}triquement
  finis.
\newblock {\em Israel J. Math.}, 118:109--124, 2000.

\bibitem[FH19]{FisHas:19}
T.~Fisher and B.~Hasselblatt.
\newblock {\em Hyperbolic Flows}.
\newblock Zurich Lectures in Advanced Mathematics. European Mathematical
  Society Publishing House, 2019.

\bibitem[Ham89]{Ham:89}
Ursula Hamenst\"{a}dt.
\newblock A new description of the {B}owen-{M}argulis measure.
\newblock {\em Ergodic Theory Dynam. Systems}, 9(3):455--464, 1989.

\bibitem[HK95]{HanKit:95}
Michael Handel and Bruce Kitchens.
\newblock Metrics and entropy for non-compact spaces.
\newblock {\em Israel J. Math.}, 91(1-3):253--271, 1995.
\newblock With an appendix by Daniel J. Rudolph.

\bibitem[Jar29]{Ja:29}
Vojtech Jarnik.
\newblock Diophantischen approximationen und hausdorffsches mass.
\newblock {\em Mat. Sb.}, 36:371--382, 1929.

\bibitem[KH95]{KatHas:95}
Anatole Katok and Boris Hasselblatt.
\newblock {\em Introduction to the modern theory of dynamical systems},
  volume~54 of {\em Encyclopedia of Mathematics and its Applications}.
\newblock Cambridge University Press, Cambridge, 1995.
\newblock With a supplementary chapter by Katok and Leonardo Mendoza.

\bibitem[Kle98]{Kle:98}
Dmitry~Y. Kleinbock.
\newblock Nondense orbits of flows on homogeneous spaces.
\newblock {\em Ergodic Theory Dynam. Systems}, 18(2):373--396, 1998.

\bibitem[KM96]{KleMar:96}
D.~Y. Kleinbock and G.~A. Margulis.
\newblock Bounded orbits of nonquasiunipotent flows on homogeneous spaces.
\newblock In {\em Sina\u{\i}'s {M}oscow {S}eminar on {D}ynamical {S}ystems},
  volume 171 of {\em Amer. Math. Soc. Transl. Ser. 2}, pages 141--172. Amer.
  Math. Soc., Providence, RI, 1996.

\bibitem[KW13]{KleWei:13}
Dmitry Kleinbock and Barak Weiss.
\newblock Modified {S}chmidt games and a conjecture of {M}argulis.
\newblock {\em J. Mod. Dyn.}, 7(3):429--460, 2013.

\bibitem[OP04]{OP:04}
Jean-Pierre Otal and Marc Peign\'{e}.
\newblock Principe variationnel et groupes kleiniens.
\newblock {\em Duke Math. J.}, 125(1):15--44, 2004.

\bibitem[PP14]{ParkPau:14}
Jouni Parkkonen and Fr\'{e}d\'{e}ric Paulin.
\newblock Skinning measures in negative curvature and equidistribution of
  equidistant submanifolds.
\newblock {\em Ergodic Theory Dynam. Systems}, 34(4):1310--1342, 2014.

\bibitem[PPS15]{PauPolSch:15}
Fr\'{e}d\'{e}ric Paulin, Mark Pollicott, and Barbara Schapira.
\newblock Equilibrium states in negative curvature.
\newblock {\em Ast\'{e}risque}, (373):viii+281, 2015.

\bibitem[PS71]{PugShu:71}
Charles Pugh and Michael Shub.
\newblock Ergodic elements of ergodic actions.
\newblock {\em Compositio Math.}, 23:115--122, 1971.

\bibitem[Riq18]{Riq:18}
Felipe Riquelme.
\newblock Ruelle's inequality in negative curvature.
\newblock {\em Discrete Contin. Dyn. Syst.}, 38(6):2809--2825, 2018.

\bibitem[RV22]{RiVe:22}
Felipe Riquelme and An\'ibal Velozo.
\newblock On the hausdorff dimension of geodesics that escape on average.
\newblock {\em In preparation}, 2022.

\bibitem[Tho11]{Tho:11}
Daniel~J. Thompson.
\newblock A thermodynamic definition of topological pressure for non-compact
  sets.
\newblock {\em Ergodic Theory Dynam. Systems}, 31(2):527--547, 2011.

\end{thebibliography}

\end{document}